\documentclass[aap,MSNbibl,dvips]{arximspdf}
\usepackage{accents}
\usepackage{graphicx}
\doi{10.1214/12-AAP860} 
\volume{23}
\issue{3}
\pubyear{2013}
\firstpage{923}
\lastpage{956}

\makeatletter

\newcommand{\rref}[1]{\fontsize{8.36pt}{10pt}\selectfont{\mbox{\ref{#1}}}}

\renewcommand{\mathring}[1]{\accentset{\circ}{#1}}

\newcommand{\shortrightarrow}{\rightarrow}

\newtheorem{theorem}{Theorem}[section]
\newtheorem{corollary}[theorem]{Corollary}
\newtheorem{proposition}[theorem]{Proposition}
\newtheorem{lemma}[theorem]{Lemma}

\newproclaim{remark}[theorem]{Remark}

\newcommand{\Iedges}{{\hat{\mathcal{I}}}}

\def\bfPu{{\mathbf P}_{u}}
\def\P{{\mathbf P}}
\def\E{{\mathbf E}}
\def\lacei{{\mathcal I}}
\def\hP{{\mathbb P}}
\def\autg{\operatorname{Aut}(G)}

\makeatother

\begin{document}
\begin{frontmatter}

\title{Random interlacements and amenability}
\runtitle{Random interlacements and amenability}

\begin{aug}
\author[A]{\fnms{Augusto} \snm{Teixeira}\corref{}\thanksref{t1}\ead[label=e1]{augusto@impa.br}}
\and
\author[B]{\fnms{Johan} \snm{Tykesson}\thanksref{t2}\ead[label=e2]{johan.tykesson@math.uu.se}}
\runauthor{A. Teixeira and J. Tykesson}
\affiliation{\'Ecole Normale Sup\'erieure and Instituto Nacional de
Matem\'atica, and~Weizmann Institute of Science}
\address[A]{\'Ecole Normale Sup\'erieure\\
D\'{e}partement de Math\'{e}matiques\\
\quad et Applications\\
Paris 75230\\
France\\
and\\
Instituto Nacional de Matem\'atica\\
Pura e Aplicada (IMPA)\\
Rio de Janeiro 22460-320\\
Brazil\\
\printead{e1}}
\address[B]{Weizmann Institute of Science\\
Faculty of Computer Science\\
\quad and Mathematics\\
Rehovot 76100\\
Israel\\
\printead{e2}} 
\end{aug}

\thankstext{t1}{Supported by a fellowship provided by AXA-Research
Fund.}

\thankstext{t2}{Supported by a post-doctoral grant of the Swedish
Research Council.}

\received{\smonth{3} \syear{2011}}
\revised{\smonth{3} \syear{2012}}

%
\begin{abstract}
We consider the model of random interlacements on transient graphs,
which was first introduced by Sznitman [\textit{Ann. of Math.} (2)
(2010) \textbf{171} 2039--2087] for the special case of ${\mathbb Z}^d$
(with $d \geq3$). In Sznitman [\textit{Ann. of Math.} (2)
(2010) \textbf{171} 2039--2087], it was shown that on ${\mathbb
Z}^d$: for any intensity $u>0$, the interlacement set is almost surely
connected. The main result of this paper says that for transient,
transitive graphs, the above property holds if and only if the graph is
amenable. In particular, we show that in nonamenable transitive graphs,
for small values of the intensity $u$ the interlacement set has
infinitely many infinite clusters. We also provide examples of
nonamenable transitive graphs, for which the interlacement set becomes
connected for large values of $u$. Finally, we establish the
monotonicity of the transition between the ``disconnected'' and the
``connected'' phases, providing the uniqueness of the critical value
$u_c$ where this transition occurs.
\end{abstract}

%

\end{frontmatter}

\section{Introduction}\label{sintro}

In~\cite{Szn07}, Sznitman introduced the model of random
interlacements on $\mathbb{Z}^d$, for $d \geq3$. Intuitively speaking,
this model describes
the local picture left by the trace of a random walk on a discrete
torus or a discrete cylinder; see~\cite{Win08} and~\cite{Szn09b}.
Moreover, recent works have shown that random interlacements can be used
to obtain a better understanding of the trace of a random walk on these
graphs; see, for instance,~\cite{Sznitman2009a,Szn09a} and
\cite{TW10}.

Later in~\cite{Tei09} the construction of random interlacements was
generalized to any transient weighted graph. This extension has
already been useful to prove results concerning the random walk on
random regular graphs; see~\cite{CTW10}. There are also strong
indications that a better understanding of the random interlacements
model in more general classes of graphs could provide interesting
results on the disconnection time of a discrete cylinder with
general basis by a random walk; see \mbox{\cite{Szn08,Szn11}} and
\cite{Win10}.

In this paper we study random interlacements on a given graph
$G=(V,E)$, with special emphasis on the relations between the behavior
of the interlacement set and geometric properties of $G$.
First, let us give an intuitive description of
the model which will be made precise in the next section. To define random
interlacements, one first considers the space $W^*$ of doubly infinite,
transient
trajectories in $G$ modulo time shift; see (\ref{eqW*}). The random
interlacements consist of a Poisson point process on the space $W^*$,
constructed in a probability space $(\Omega, \mathcal{A},
\mathbb{P}_u)$. Usually one is interested in understanding the trace
left by this soup of random walk trajectories, which is a random subset
$\mathcal{I}$ of $V$. The parameter $u \in
\mathbb{R}_+$ is a multiplicative factor of the intensity; therefore
the bigger the value of $u$, the more trajectories enter the picture.
Although we postpone the precise definition of random interlacements
to the next section, we mention here that $\mathcal{I}$, under the law
$\mathbb{P}_u$, is
characterized as the unique random subset of $V$ such that
%
\begin{equation}
\label{echarIu} \mathbb{P}_u[\mathcal{I} \cap K = \varnothing] =
\exp\bigl\{-u \operatorname{cap}(K)\bigr\}\qquad \mbox{for every finite $K \subset
V$};
\end{equation}
see (1.1) of~\cite{Tei09}. In the equation above,
$\operatorname{cap}(K)$ stands for the capacity of $K$ defined in
(\ref{eqcapeq}).

In the case of $\mathbb{Z}^d$, it was already shown that
$\mathcal{I}$, regarded as a random subset of~$\mathbb{Z}^d$, defines
an ergodic percolation process; see, for instance,~\cite{Szn07},
Theorem 2.1. However, due to its long-range dependence, $\mathcal{I}$
behaves very differently from the usual independent site percolation
on $\mathbb{Z}^d$. Let us stress here one such difference.
Although the marginal density ${\mathbb P}_u[x\in\mathcal{I}]$
converges to zero as
we drive the intensity parameter $u$ to zero, the set ${\mathcal I}$ is
${\mathbb P}_u$-almost surely given by a
single infinite connected subset of $\mathbb{Z}^d$; see~\cite{Szn07},
Corollary~2.3. In particular there is no phase transition for the
connectivity of $\mathcal{I}$ as one varies $u$. For this reason, in the
case of~$\mathbb{Z}^d$, more attention has been given to the complement
of $\mathcal{I}$, the so called vacant set denoted by ${\mathcal V}$,
which exhibits a
nontrivial phase transition for every $d \geq3$; see, for instance,
\cite{Szn07} and
\cite{SS09}.

In this work, we address similar questions for more general graphs
$G$. We first show that
%
\begin{equation}
\label{eergodic}
\begin{tabular}{p{323pt}}
if $G$ is a transient, vertex-transitive graph, then the
interlacement set $\mathcal{I}$ under the law $\mathbb{P}_u$
is an ergodic subset of $V$;
\end{tabular}\hspace*{-27pt}
\end{equation}
see Proposition~\ref{t01law}. In particular, any
automorphism-invariant event on
$\{0,1\}^{V}$ has probability either $0$ or $1$.

Given this fact, it is natural to ask what special property of $\mathbb
{Z}^d$ is responsible for the connectivity of the interlacement set at
all intensities $u$. In this paper we show that the answer to this
question is related to the amenability of $\mathbb{Z}^d$, as shown in
the two results that we now describe. The main result of this paper is
the following:
%
\begin{equation}
\label{emain}
\begin{tabular}{p{323pt}}
if $G$ is a vertex-transitive nonamenable graph,
then for $u > 0$ small enough,
the interlacement set $\mathcal{I}$ is $\mathbb{P}_u$-almost
surely disconnected;
\end{tabular}\hspace*{-27pt}
\end{equation}
see Theorem~\ref{tdisconn}. This shows a clear difference from the
$\mathbb{Z}^d$ case.

This result motivates the definition of the critical parameter
%
\begin{equation}
\label{ucdef} u_c=u_c(G)=\inf\{u\dvtx  \lacei\mbox{ is
connected } {\mathbb P}_u\mbox{-a.s.}\}.
\end{equation}
See also Corollary~\ref{coneuc}. A trivial consequence of (\ref
{emain}) is that the critical threshold
$u_c$ is positive for vertex-transitive, nonamenable graphs. In fact,
Theorem~\ref{tdisconn} gives a lower bound for $u_c$ in terms of the
spectral radius of $G$; see definition in (\ref{spectraldef}).

The counterpart for Theorem~\ref{tdisconn} is given in Theorem \ref
{tamenable}, which gives the following generalization of Corollary 2.3
of~\cite{Szn07}:
%
\begin{equation}
\label{eamenstat}
\begin{tabular}{p{323pt}}
if $G$ is a vertex-transitive amenable transient graph,
then for all $u>0$, the set ${\mathcal I}$ is ${\mathbb
P}_u$-almost surely connected.
\end{tabular}\hspace*{-27pt}
\end{equation}
In particular, (\ref{eamenstat}) implies that for transitive amenable
graphs $G$, $u_c(G)=0$.

We also note that it is not always the case that $u_c$ is finite. In
fact according to Remark~\ref{rtreeinf}, in an infinite regular tree
$\mathbb{T}^d$ for $d \geq3$, the interlacement set is almost surely
disconnected
for every value of $u > 0$, so that $u_c(\mathbb{T}^d)=\infty$.
However, in the same spirit as for Bernoulli
percolation (see, e.g.,~\cite{GN90}), we can prove that
%
\begin{equation}
\label{etreez}
\begin{tabular}{p{323pt}}
for random interlacements on $\mathbb{T}^d
\times \mathbb{Z}^{d'}$, where $d \geq3$ and $d' \geq1$,
we have $0 < u_c (\mathbb{T}^d \times
\mathbb{Z}^{d'}) < \infty$;
\end{tabular}\hspace*{-27pt}
\end{equation}
see Proposition~\ref{pucfinite}.
Thus each of the three cases $u_c(G)=0$, $u_c(G)\in(0,\infty)$ and
$u_c(G)=\infty$ can happen.
It is also of interest to study what happens regarding the connectivity
of ${\mathcal I}$ above $u_c$ (note that the connectivity of $\mathcal
{I}$ is not a monotone property). In Theorem~\ref{tmonotone}, we show
that
%
\begin{equation}
\label{euniqmonotone}
\begin{tabular}{p{323pt}}
if $G$ is a vertex-transitive and transient graph,
then for any $u>u_c$,
the set ${\mathcal I}$ is ${\mathbb P}_u$-almost surely
connected.
\end{tabular}\hspace*{-27pt}
\end{equation}
In other words, there is monotonicity in the uniqueness transition.

In addition to the critical value $u_c$ defined in (\ref{ucdef}),
there is another critical value of interest, which was defined in \cite
{Szn07}, (0.13). Recall that ${\mathcal V}$ stands for the complement
of ${\mathcal I}$. For $x \in V$, let
%
\begin{eqnarray}
\label{eetadef} \eta_x(u)&=&{\mathbb P}_u[x\mbox{ belongs
to an unbounded connected component of } {\mathcal V}],
\\
%
\label{eustardef} u_* (G)&=&\inf\bigl\{u\ge0,\eta_x(u)=0\bigr\}.
\end{eqnarray}
In Corollary 3.2 of~\cite{Tei09}, the critical intensity $u_*(G)$ was
shown to be independent of the choice of the base point $x$.

The critical intensity of $\mathbb{Z}^d$ was first studied in \cite
{Szn07}, where it was shown that $u_*({\mathbb Z}^d)<\infty$ for $d\ge
3$, and $u_*({\mathbb Z}^d)>0$ for $d\ge7$. Later in~\cite{SS09},
$u_*({\mathbb Z}^d)>0$ was established for any $d\ge3$. In \cite
{Tei09}, the nondegeneracy of $u_*$ was established for some families
of graphs.

In Proposition~\ref{puthresbound} we join results from~\cite{Tei09,BLPS97} and (\ref{eergodic}) to show that
%
\begin{equation}
\label{eustar} 0 < u_*(G) < \infty\qquad\mbox{whenever $G$ is a nonamenable Cayley
graph.}
\end{equation}

We now give a brief overview of the proof of the main result (\ref
{emain}) of this paper. We rely on the following
domination result:
%
\begin{equation}
\label{edomin}
\begin{tabular}{p{323pt}}
there is a coupling between $\mathbb{P}_u$ and
the law of a certain
branching random walk on $G$, such that the cluster of $\mathcal {I}$
containing a given point $x$ is almostsurely contained in
the trace left by the branching random walk from~$x$.
\end{tabular}\hspace*{-27pt}
\end{equation}
Given this coupling, the proof of (\ref{emain}) is reduced to a
calculation on the heat-kernel of the branching random walk.

The domination stated in (\ref{edomin}) is given in two steps. First,
in Proposition~\ref{pdomin} we give the domination of random
interlacements by the ``frog model'' or ``$A + B \to2A$ model,''
studied in~\cite{AMP02} and~\cite{RS04}. Then
Proposition~\ref{pdomin2} establishes the domination of the mentioned
frog model by a specific branching random walk on~$G$.

We now discuss the resemblance between our main results and a
conjecture in Bernoulli percolation. As we explained above, our results
characterize the amenability of transitive graphs in terms of the
existence of infinitely many infinite clusters in the interlacements
set, for some parameter $u > 0$. Our investigation was partially
inspired by a similar, still open question for Bernoulli percolation.
More precisely, Benjamini and Schramm conjectured that a transitive
graph is nonamenable if and only if there is some $p \in[0,1]$ for
which Bernoulli site percolation with parameter $p$ contains an
infinite number of infinite clusters; see Conjecture 6 in \cite
{Benjamini1996}. This conjecture has been resolved in some cases:
Benjamini and Schramm~\cite{BenSchr2001} solved it in the planar case
and Pak and Smirnova-Nagnibeda~\cite{PakSmirnova2000} established the
nonuniqueness phase for certain classes of nonamenable Cayley graphs.
Theorem 4 of~\cite{Benjamini1996} gives a sufficient condition for
having a nonuniqueness phase, and that theorem is the main inspiration
for our result concerning disconnectedness of random interlacements. In
\cite{USF2001}, Benjamini et al. characterized the amenability property
of Cayley graphs in terms of the connectedness of the union of the
wired spanning forest and Bernoulli percolation.

The rest of the paper is organized as follows. The notation and the main
properties of random interlacements we need are given in
Section~\ref{snotation}, where also the $\{0,1\}$-law is provided. We
then prove that the interlacement set is always connected on amenable
graphs; see Section~\ref{samenablecase}. The Propositions \ref
{pdomin} and
\ref{pdomin2} concerning domination of random interlacements by the
frog model and the branching random walk are given, respectively, in
Sections~\ref{ssleepart} and~\ref{sbranching}.
The proof of Theorem~\ref{tdisconn} concerning disconnectedness of the
interlacement set for small $u$ is given in Section~\ref
{sdisconnsection}. In the same section, we also obtain the required
heat-kernel estimates on the mentioned branching random walk. In
Section~\ref{sconnectedsection}, we prove that $u_c$ is finite for the
graph $\mathbb{T}^d \times\mathbb{Z}^{d'}$ and deal with the
monotonicity of the uniqueness transition. Finally, in Section \ref
{sulowerbound}, we obtain bounds on $u_*$ for nonamenable Cayley
graphs, establishing~(\ref{eustar}).

In this paper we use the following convention on constants, which are
denoted by $c$ and $c'$. These constants can change from line to line
and solely depend on the given graph $G$ which is fixed within proofs.
Further dependence of constants are indicated by subscripts, that is,
$c_\varepsilon$ is a constant depending on $\varepsilon$ and potentially on
the graph $G$.

\section{Notation and definitions}
\label{snotation} We let $G=(V,E)$ be an infinite connected graph,
which we always assume to have finite geometry, that is, every vertex
in $V$ belongs only to a finite number of edges. For $x,y \in V$, we
write $x \leftrightarrow y$ if they are neighbors in $G$. For $x\in V$,
let $d_x$ be the degree of $x$. Let $\Delta=\Delta(G)=\sup_{x\in V}
d_x$, which sometimes will be assumed to be finite. We recall the
definitions of vertex-transitivity and amenability. A bijection
$f\dvtx  V\to V$ such that $\{f(x),f(y)\}\in E$ if and only if
$\{x,y\}\in E$ is said to be a graph automorphism of $G$. We denote the
set of automorphisms of $G$ by $\operatorname{Aut}(G)$ and say that
$G$ is
transitive if for any $x,y\in V$, there is a graph automorphism $f$
such that $f(x)=y$. For $K\subset V$, we define the \textit{outer vertex
boundary} of $K$ as $\partial_V K=\{y\notin K\dvtx  \exists x\in K
,d(x,y)=1\}$, where $d$ is the usual graph distance induced in $G$. We
write $\overline K$ for the set $K \cup\partial_V K$. {We also define the
\textit{outer edge boundary} of $K$ as $\partial_E K=\{\{x,y\}\in
E\dvtx  x\in K, y\notin K\}$.}

The \textit{vertex isoperimetric constant} of the graph $G$ is defined as
%
\begin{equation}
\label{ekvdef}\kappa_V=\kappa_V(G)=\inf \biggl\{
\frac
{|\partial_V K|}{|K|}\dvtx  |K|<\infty \biggr\}.
\end{equation}
Similarly, the \textit{edge isoperimetric constant} is defined to be
%
\begin{equation}
\label{ekedef}\kappa_E=\kappa_E(G)=\inf \biggl\{
\frac
{|\partial_E K|}{|K|}\dvtx  |K|<\infty \biggr\}.
\end{equation}
A graph $G$ is said to be nonamenable if $\kappa_V(G)>0$, otherwise it
is said to be amenable. In Theorems~\ref{tamenable} and \ref
{tdisconn}, where we use the notion of amenability, we also assume the
graph to have bounded degree; therefore we could have defined
amenability in terms of $\kappa_E$ instead.

We let $P_x$ stand for the law of a simple random walk on $G$ starting
at $x$, while $(X_n)_{n \geq0}$ stands for the canonical projections
on $V$.
For vertices $x,y\in V$, let
$p^{(n)}(x,y)$ be the probability that a simple random walk
started at $x$ is at $y$ after $n$ steps, that is, $P_x[X_n = y]$. Let
%
\begin{equation}
\label{spectraldef}\rho=\rho(G):=\limsup_{n\to
\infty} \bigl(p^{(n)}(x,y)
\bigr)^{1/n}.
\end{equation}
The quantity $\rho$ is called the
\textit{spectral radius} of $G$, and is independent of the choices
of $x$ and $y$. If $G$ is nonamenable and has bounded degree, then
$\rho<1$; see
\cite{D84}, Theorem 2.3. We also define the Green's function of the
simple random walk on $G$ as $g(x,y) = \sum_{n \geq0} p^{(n)}(x,y)$,
for $x, y \in V$.

The space $W_+$ stands for the set of infinite trajectories that spend
only a finite time in finite sets
%
\begin{eqnarray}
W_+ &=& \bigl\{ \gamma\dvtx  \mathbb{N} \rightarrow V; \gamma(n) \leftrightarrow
\gamma(n+1) \mbox{ for each } n \geq 0 \mbox{ and }
\nonumber\\[-8pt]\\[-8pt]
&&\hspace*{64.1pt}\bigl\{n; \gamma(n) = y\bigr\} \mbox{ is finite for all } y \in V \bigr\}.
\nonumber
\end{eqnarray}
We endow $W_+$ with the $\sigma$-algebra $\mathcal{W}_+$ generated by
the canonical coordinate maps~$X_n$.

We further consider the space of doubly infinite trajectories that
spend only a finite time in finite subsets of $V$
%
\begin{eqnarray}
W &=& \bigl\{ \gamma\dvtx \mathbb{Z}
\rightarrow V; \gamma(n) \leftrightarrow \gamma(n+1)\mbox{ for each } n \in
\mathbb{Z}\mbox{ and}
\nonumber\\[-8pt]\\[-8pt]
&&\hspace*{66pt}\bigl\{n; \gamma(i) = y\bigr\} \mbox{ is finite for all } y \in V \bigr\}.
\nonumber
\end{eqnarray}

On the space $W$, for $k \in\mathbb{Z}$, we introduce the shift
operator $\theta_k\dvtx  W \to W$ which sends a trajectory $w$ to $w'$ such
that $w'(\cdot) = w(\cdot- k)$. We also consider the space $W^*$ of
trajectories in $W$ modulo time shift
%
\begin{equation}
\label{eqW*} W^* = W / \sim\qquad \mbox{where } w \sim w' \iff w =
\theta_k\bigl(w'\bigr) \qquad\mbox{for some } k \in
\mathbb{Z}\hspace*{-10pt}
\end{equation}
and denote with $\pi^*$ the canonical projection from $W$ to $W^*$. The
map $\pi^*$ induces a $\sigma$-algebra in $W^*$ given by $\mathcal{W}^*
= \{A \subset W^*; (\pi^*)^{-1}(A) \in\mathcal{W}\}$, which is the
largest $\sigma$-algebra on $W^*$ for which $(W,\mathcal{W})
\stackrel
{\pi^*}{\rightarrow} (W^*,\mathcal{W}^*)$ is measurable.

For any finite set $K \subset V$, we define for a trajectory $w \in W$
the entrance time of $K$ as
%
\begin{equation}
\label{eqhitK} H_K(w) = \inf\bigl\{k \in\mathbb{Z}; w(k) \in K\bigr
\}.
\end{equation}
For a trajectory $w \in W_+$, the hitting time of $K$ is defined as
%
\begin{equation}
\tilde{H}_K(w) = \inf\bigl\{k \ge1; w(k) \in K\bigr\}.
\end{equation}

Considering still a finite $K \subset V$, we define the \textit
{equilibrium measure} $e_K$ by
%
\begin{equation}
\label{eqequilibmeasure} e_{K}(x) = 1_{\{x \in K\}}P_x[
\tilde{H}_K = \infty] \cdot{d_x}
\end{equation}
and the capacity
%
\begin{equation}
\label{eqcapeq} \operatorname{cap}(K) = \sum_{x \in K}
e_K(x).
\end{equation}
In addition, we introduce the measure
%
\begin{equation}
\label{epekdef} P_{e_K}=\sum_{x\in V}e_K(x)P_{x}.
\end{equation}

Given a finite set $K \subset V$, write $W_K$ for the space of
trajectories in $W$ that enter the set $K$, and denote with $W_K^*$ the
image of $W_K$ under $\pi^*$.

The set of point measures on which one canonically defines random
interlacements is given by
%
\begin{equation}
\label{eqomega} \Omega= \biggl\{\omega= \sum_{i\geq 1}
\delta_{w^*_i}; w^*_i \in W^* \mbox{ and } \omega
\bigl(W^*_K\bigr) < \infty\mbox{, for every finite } K \subset V
\biggr\},\hspace*{-35pt}
\end{equation}
endowed with the $\sigma$-algebra $\mathcal{A}$ generated by the
evaluation maps $\omega\mapsto\omega(D)$ for $D \in\mathcal{W}^*
\otimes\mathcal{B}(\mathcal{R}_+)$.

The following theorem resembles Theorem 1.1 in~\cite{Szn07}. It was
established in Theorem 2.1 of~\cite{Tei09}, providing the existence of
the intensity measure used to construct the random interlacements process.
%
\begin{theorem}
\label{thexistnu}
There exists a unique $\sigma$-finite measure $\nu$ on $(W^*,\mathcal
{W}^*)$ satisfying, for each finite set $K \subset V$,
%
\begin{equation}
\label{eqnuQ} 1_{W^*_K} \cdot\nu= \pi^* \circ Q_K,
\end{equation}
where the finite measure $Q_K$ on $W_K$ is determined by the following.
Given $A$ and $B$ in $\mathcal{W}_+$ and a point $x \in V$,
%
\begin{eqnarray}
\label{eqQK}
&&
Q_K\bigl[(X_{-n})_{n\geq 0} \in A,
X_0 = x, (X_n)_{n\geq 0} \in B\bigr] \nonumber\\[-8pt]\\[-8pt]
&&\qquad=
P_x[A|{\tilde{H}}_K = \infty] e_K(x)
P_x[B].\nonumber
\end{eqnarray}
\end{theorem}

We are now ready to define the random interlacements. Consider on
$\Omega$ the law $\mathbb{P}_u$ of a Poisson point process with
intensity measure given by $u \nu(dw^*)$; for a reference on the
construction of Poisson point processes, see~\cite{R08}, Proposition
3.6. We define the \textit{interlacement} and the \textit{vacant set}
at level $u$, respectively, as
%
\begin{equation}
\label{eqinterlace} \mathcal{I} (\omega) = \biggl\{ \bigcup
_{w^* \in\operatorname
{supp}(\omega)} \operatorname{Range}\bigl(w^*\bigr) \biggr\}
\end{equation}
and
\begin{equation}
\label{eqvacant} \mathcal{V}(\omega) = V \setminus\mathcal{I} (\omega)
\end{equation}
for $\omega= \sum_{i\geq 0} \delta_{w^*_i}$ in $\Omega$.\vspace*{1pt}

Observe that with (\ref{echarIu}), (\ref{eqequilibmeasure}) and
(\ref{eqcapeq}), we have
%
\begin{equation}
\label{evuprobx} {\mathbb P}_u[x\in{\mathcal V}]=\exp\bigl\{-u
d_x P_x[\tilde {H}_x=\infty ]\bigr\}.
\end{equation}

For finite $K\subset V$ and $\omega=\sum_{i\ge0}\delta_{w_i^*}\in
\Omega$ we define
%
\begin{equation}
\label{emukdefin} \mu_K(\omega)=\sum_{i\ge0}
\delta_{(w_i^*)^{K,+}} 1\bigl\{w_i^*\in W_K^*\bigr\}
\end{equation}
were for $w^*\in W^*_K$, $(w_i^*)^{K,+}$ is defined to be the
trajectory in $W_+$ which follows $w^*$ from the first time it hits $K$.
Proposition 1.3, equation (1.45) in~\cite{Szn07} says that
%
\begin{equation}
\label{emukthm}
\begin{tabular}{p{323pt}}
under ${\mathbb P}_u$, $\mu_K$ is
a Poisson point process on $W_+$
with intensity measure $u P_{e_K}(dw)$.
\end{tabular}\hspace*{-27pt}
\end{equation}
Actually, in~\cite{Szn07}, (\ref{emukthm}) is stated only for the
${\mathbb Z}^d$ case, but the proof does not use any special properties
of ${\mathbb Z}^d$ and is thus also valid here.
%
\begin{remark}
Recall the definition of $H_K$ in (\ref{eqhitK}) and note
that for each $w^* \in W^*_K$, there is a unique $w \in W_K$ such that
$\pi^*(w) = w^*$ and $H_K(w) = 0$, we call this $s_K(w^*)$. Therefore,
since by (\ref{eqQK}) the measure $Q_K$ is supported on $\{H_K = 0\}
$, we conclude that for any $A \in\mathcal{W}$ contained in $W_K$,
%
\begin{eqnarray}
\label{eqsuppQK}
&&
Q_K \bigl[\bigl(\pi^*\bigr)^{(-1)} \bigl(
\pi^*(A)\bigr)\bigr]\nonumber\\
&&\qquad= Q_K \bigl[ \bigl\{ w \in W;
\theta_l (w) \in A \mbox{ for some $l \in\mathbb{Z}$}\bigr\} \bigr]
\nonumber\\[-8pt]\\[-8pt]
&&\qquad= Q_K \bigl[ \bigl\{ w \in W; H_K(w) = 0 \mbox{ and }
\theta_l (w) \in A \mbox{ for some $l \in\mathbb{Z}$}\bigr\} \bigr]
\nonumber
\\
&&\qquad= Q_K \bigl[s_K\bigl(\pi^*(A)\bigr)\bigr].\nonumber
\end{eqnarray}
\end{remark}

For a measure $\gamma$ on $W^*$ and $g\in\autg$, let $\gamma^g$ denote
the image of $\gamma$ under the mapping $w^*\to g(w^*)$. In Theorem 1.1
of~\cite{Szn07} it was shown that on ${\mathbb Z}^d$, the measure $\nu$
is invariant under translations. We now show a similar statement for
vertex-transitive graphs.
%
\begin{proposition}\label{pnuinvariance}
For $g\in\autg$,
%
\begin{equation}
\label{enuinvariance} \nu^g=\nu.
\end{equation}
\end{proposition}
\begin{pf}
We follow the proof of Theorem 1.1 in~\cite{Szn07}. For $g\in\autg$
and $K\subset V$ finite, $w^*\to g(w^*)$ maps $W_K^*$ one-to-one onto
$W^*_{g(K)}$. Also note that $s_K(g(w^*))=g(s_{g^{-1}(K)}(w^*))$ for
every $w^*\in W^*_{g^{-1}(K)}$. Then, for $C\in{\mathcal W}$, we get that
%
\begin{eqnarray}
\label{eskg} s_K\circ\bigl(1_{W^*_K} \cdot\nu^g
\bigr) (C) & = & s_K\circ \bigl((1_{W^*_{g^{-1}(K)}} \cdot\nu)^g
\bigr) (C)
\nonumber\\
&=&s_{g^{-1}(K)}\circ(1_{W^*_{g^{-1}(K)}} \cdot\nu) \bigl(g(w)\in C \bigr)
\\
&=& Q_{g^{-1}(K)}\bigl[g(w)\in C\bigr].\nonumber
\end{eqnarray}
Now fix $A,B\in{\mathcal W}_+$ and $y\in V$. Let $C$ be the event $\{
(X_{-n})_{n\geq 0} \in A, X_0 = y, (X_n)_{n\geq 0} \in B\}$.
Using (\ref{eqQK}), (\ref{eqnuQ}) and (\ref{eskg}), we obtain that
%
\begin{eqnarray}
&&s_K \circ\bigl(1_{W^*}\nu^g\bigr) (C)
\nonumber
\\
&&\qquad =Q_{g^{-1}(K)} \bigl[(X_{-n})_{n\geq 0} \in
g^{-1}(A), X_0 = g^{-1}(y), (X_n)_{n\geq 0}
\in g^{-1}(B) \bigr]
\nonumber\\[-8pt]\\[-8pt]
&&\qquad=P_{g^{-1}(y)}^{g^{-1}(K)}\bigl[g^{-1}(A)\bigr]e_{g^{-1}(K)}
\bigl(g^{-1}(y)\bigr)P_{g^{-1}(y)}\bigl[g^{-1}(B)
\bigr]=Q_K[C]
\nonumber
\\
&&\qquad=s_K\circ(1_{W^*_K}\nu) (C).\nonumber
\end{eqnarray}
This clearly implies that $\nu^g=\nu$.
\end{pf}
For $g\in\autg$, we define $\tau_g\dvtx  \Omega\to\Omega$, which maps
$\omega=\sum_{i\ge0}\delta_{w_i^*}\in\Omega$ to
%
\begin{equation}
\label{eshift} \tau_g\omega=\sum_{i\ge0}
\delta_{g(w_i^*)}\in\Omega.
\end{equation}
A consequence of (\ref{enuinvariance}) is that
%
\begin{equation}
\label{epuinvariance} {\mathbb P}_u\mbox{ is invariant under } (
\tau_g)_{g\in\autg}.
\end{equation}

More than that, we can prove a $\{0,1\}$-law. First, we need to
introduce some additional notation. A given $g\in\autg$ naturally
induces a mapping from $\{0,1\}^V$ to itself which we shall denote by
$t_g$. Denote by $Q_u$ the law on $\{0,1\}^V$ of $(1\{x\in\mathcal
{I}\}
)_{x\in V}$ under ${\mathbb P}_u$. We write ${\mathcal Y}$ for the
$\sigma$-algebra generated by the coordinate maps $Y_z$, $z\in V$ on
$\{
0,1\}^V$.
%
\begin{proposition}\label{t01law}
If $G$ is transient and transitive, then for any $u\ge0$, $(t_g)_{g\in
\autg}$ is a measure preserving flow on $(\{0,1\}^V,{\mathcal Y},Q_u)$
which is ergodic. In particular, for any $u\ge0$ and any $A\in
{\mathcal Y}$ which is invariant under $\operatorname{Aut}(G)$,
%
\begin{equation}
\label{e01law}Q_u[A]\in\{0,1\}.
\end{equation}
\end{proposition}

The proof of this proposition is very similar to the proof of
Theorem 2.1 in~\cite{Szn07}, but for the reader's convenience we
present the full proof with the necessary modifications in the
\hyperref[app]{Appendix}.

\section{The amenable case}\label{samenablecase}

In this section we prove that if $G$ is an a\-me\-na\-ble transitive
graph, then $\mathcal{I}$ is ${\mathbb P}_u$-almost surely connected
for every $u > 0$. This is not a direct consequence of the results in
\cite{BK89} as we explain in Remark~\ref{runique} below.
%
\begin{remark}
\label{runique}
In~\cite{BK89}, Burton and Keane developed a quite general technique to
prove uniqueness of the infinite components induced by a random subset
of edges $\mathcal{U} \subset E$ of an amenable graph $G = (V,E)$; see
also Theorem~12.2 in~\cite{HJ06}. Their result applies under\vadjust{\goodbreak} the rather
general conditions that $\mathcal{U}$ is translation invariant and
satisfies the so called \textit{finite energy property}; that is, for
every edge $e \in E$ and any $U \subset E \setminus\{e\}$,
%
\begin{equation}
\label{efiniteen} 0 < \mathcal{P}\bigl[e \in\mathcal{U} | \mathcal{U} \setminus
\{e\} = U\bigr] < 1
\end{equation}
for some version of this conditional probability; see Definition 12.1
in~\cite{HJ06}. We note here that this condition does not hold in
general for the interlacement set $\mathcal{I}$ as observed in
Remark 2.2(3) of~\cite{Szn07}.
\end{remark}

Instead of considering $\mathcal{I}$, it will be convenient to consider
the set of edges traversed by any of the trajectories in the Poisson
point process $\omega$. More precisely, consider the set
%
\begin{eqnarray}
\label{ebarI}
\Iedges(\omega) &=& \bigl\{ e \in E; e =
\bigl\{w(k),w(k+1)\bigr\},\nonumber\\[-8pt]\\[-8pt]
&&\hspace*{5.1pt}\mbox{for some $k \in\mathbb{Z}$}
\mbox{ and some $w$ such that $\pi^*(w) \in
\operatorname {supp}(\omega )$} \bigr\}.\nonumber
\end{eqnarray}

Observe, as in~\cite{Szn07} below (2.18), that the connected components
of $V$ induced by the edges in $\Iedges$ are either isolated points or
infinite sub-components of $\mathcal{I}$.
We introduce some notation as in~\cite{Szn07}. Denote by $\tilde{\psi
}\dvtx  \Omega\to\{0,1\}^{E}$ the map $\tilde{\psi}(\omega)=(1\{e\in
\Iedges(\omega)\})_{e\in E}$. Let\vspace*{1pt} $\tilde{Q}_u$ be the probability
measure on $(\{0,1\}^{E},\tilde{\mathcal Y})$ defined as the image of
${\mathbb P}_u$ under $\tilde{\psi}$. Here $\tilde{{\mathcal Y}}$
stands for the canonical $\sigma$-algebra on $\{0,1\}^E$. Let $\tilde
{t}_g$, $g\in\autg$ be the canonical shift on $\{0,1\}^E$. Then
$\tilde
{t}_g\circ\tilde{\psi}=\tilde{\psi}\circ\tau_g$, $g\in\autg$, and
%
\begin{eqnarray}
\label{eedgeergod}
\begin{tabular}{p{323pt}}
for $u\ge0$, $(\tilde{t}_g)_{g\in\autg}$
is a measure preserving flow on
$(\{0,1\}^E,\tilde{\mathcal{Y}},\tilde{Q}_u)$
which is ergodic.
\end{tabular}\hspace*{-27pt}
\end{eqnarray}
To prove (\ref{eedgeergod}), one proceeds in the same way as for
Proposition~\ref{t01law}; see the \hyperref[app]{Appendix}.

According to Remark~\ref{runique}, we need to adapt the technique of
\cite{BK89} in order to apply it to the sets $\mathcal{I}$ or
$\Iedges$.
%
\begin{proposition}\label{pnumberofclusters}
Let $G$ be a transient and transitive graph and fix $u > 0$. The number
of infinite connected components induced by the set $\hat{\mathcal{I}}$
as defined in (\ref{ebarI}) is $\mathbb{P}_u$-almost surely a
constant, which is either $1$ or $\infty$.
\end{proposition}
\begin{pf}
The proof of Corollary 2.3 of~\cite{Szn07} would go through in the
current setting with only minor modifications, but here we proceed in a
slightly different manner.
We denote by $N$ the number of connected components induced by $\hat
{\mathcal{I}}(\omega)$. The fact that $N$ is almost surely a constant
follows from ergodicity; see (\ref{eedgeergod}). Moreover, $\Iedges$
induces at least one unbounded connected component a.s. Thus it remains
to show that for any $2\le k<\infty$, ${\mathbb P}_u[N=k]=0$. Assume
for contradiction that for some $2\le k <\infty$ we have ${\mathbb
P}_u[N=k]=1$. Fix $o \in V$, under the above assumption, we can pick
$L<\infty$ so large that the event
\[
A=\bigl\{ \hat{\mathcal{I}}(\omega) \mbox{ has $k$ components, all intersecting
$B(o,L)$}\bigr\}\vadjust{\goodbreak}
\]
has positive probability. Denote by $W_1^*$ the set of trajectories in
$W^*$ that visit \textit{every} vertex in $B(o,L)$. Decompose $\omega$
into $\omega^1=1_{W_1^*}\omega$ and $\omega^2=1_{W^* \setminus W_1^*}
\cdot\omega$. Under ${\mathbb P}_u$, $\omega^1$ and $\omega^2$ are
independent Poisson point processes with intensity measures $u
1_{W_1^*}\,d\nu$ and $u 1_{W^* \setminus W_1^*} \cdot d\nu$,
respectively. Using only the transience and connectedness of $G$, it
readily follows that the mass of $1_{W_1^*}\,d\nu$ is
strictly\vspace*{1pt}
positive. Since each trajectory in $W_1^*$ intersects every trajectory
of $W_{B(o,L)}^*$, it follows that if $\omega^1(W^*)>0$, only one
connected component induced by $\hat{\mathcal{I}}(\omega)$ intersects
$B(o,L)$. Thus we have the inclusion
%
\begin{equation}
\label{eaincl} A\subset\bigl\{\omega^1\bigl(W^*\bigr)=0\bigr\}.
\end{equation}
From (\ref{eaincl}) it follows that $A\subset\tilde{A}$ where
%
\begin{equation}
\tilde{A}= \bigl\{ \hat{\mathcal{I}}\bigl(\omega^2\bigr) \mbox{ has
$k$ components, all intersecting $B(o,L)$}\bigr\}
\end{equation}
and since $A$ has positive probability the same holds for $\tilde{A}$.
Since the events $\tilde{A}$ and $\{\omega^1(W^*)>0\}$ are defined in
terms of independent Poisson point processes, they are independent. Therefore,
%
\begin{equation}
\label{eaeq2} {\mathbb P}_u\bigl[\tilde{A} \cap\bigl\{
\omega^1\bigl(W^*\bigr)>0\bigr\}\bigr]>0.
\end{equation}
But on $\tilde{A} \cap\{\omega^1(W^*)>0\}$, we have $N=1$.
Thus (\ref{eaeq2}) contradicts our assumption that $N=k$ a.s. for
some $2\le k
<\infty$, completing the proof of the proposition.
\end{pf}

In the next theorem, we rule out the possibility of having infinitely
many components in $\Iedges$ if the underlying graph is transitive and
amenable.
%
\begin{theorem}\label{tamenable}
Let $G$ be a transient, transitive and amenable graph. Then
%
\begin{equation}
\mathbb{P}_u [\mathcal{I} \mbox{ is connected}] =
\mathbb{P}_u [\Iedges\mbox{ is connected}] = 1\qquad \mbox{for every $u
> 0$.}
\end{equation}
\end{theorem}
\begin{pf}
Throughout this proof, we use the convention that Range$(w)$ is the set
of edges that are traversed by the trajectory $w$ where $w\in W$ or
$w\in W^*$. Using Proposition~\ref{pnumberofclusters}, we conclude
that for every $u > 0$, the number of infinite clusters in $\Iedges$ is
$\mathbb{P}_u$-a.s. a constant which is either $1$ or $\infty$. Since
$\Iedges$ has no finite components, all we need to do is to rule out
the case that it has infinitely many infinite clusters. For this let us
suppose by contradiction the contrary. Under this assumption, we are
going to construct, with positive probability, a trifurcation point for
$\Iedges$ as defined below.

The definition of trifurcation point was first introduced in \cite
{BK89}. We say that a given point $y$ is a trifurcation point for the
configuration $\mathcal{E} \in\{0,1\}^{E}$ if it belongs to an
infinite connected component induced by $\{\mathcal{E} = 1\}$ but the
removal of $y$ would split its cluster into three distinct infinite
connected components.

Using the assumption that $\mathbb{P}_u$-a.s. there are infinitely many
infinite connected components in $\Iedges$, we conclude that there
exist a finite, connected set $K \subseteq V$ which intersects three
distinct infinite clusters of $\Iedges$ with positive probability. This
implies that
%
\begin{eqnarray}
\label{ew123}
&&
\mathbb{P}_u \bigl[ \mbox{there are trajectories
$w_1, w_2, w_3 \in\operatorname{supp}(\omega) \cap
W^*_K$}\nonumber\\[-8pt]\\[-8pt]
&&\hspace*{18.5pt}\mbox{whose ranges belong to distinct components of $\Iedges$}
\bigr]>0.\nonumber
\end{eqnarray}

It is useful now to decompose the point measure $\omega$ into $\omega_1
= \mathbf{1}_{W^*_K} \cdot\omega$ and $\omega_2 = \mathbf{1}_{W^*
\setminus W^*_K} \cdot\omega$, which are independent under $\mathbb
{P}_u$; see the notation below (\ref{eqhitK}). We denote the laws of
$\omega_1$ and $\omega_2$, respectively, by $\mathbb{P}_1$ and
$\mathbb
{P}_2$. We say that a given set $\mathcal{U} \subseteq E$ is
\textit{good} if
%
\begin{eqnarray}
\label{eUgood}
&&
\mathbb{P}_1 \bigl[\mbox{there are $w_1$,
$w_2$, $w_3$ in $\operatorname{supp}(
\omega_1)$}\nonumber\\[-8pt]\\[-8pt]
&&\hspace*{15.5pt} \mbox{whose images are not connected through $\mathcal{U}$}
\bigr]>0.\nonumber
\end{eqnarray}
By (\ref{ew123}) and Fubini's theorem, we conclude that with positive
$\mathbb{P}_2$-probability, $\hat{\mathcal{I}}(\omega_2)$ is good.

From now on we fix a good set $\mathcal{U}$. Our aim is to show that
with positive $\mathbb{P}_1$-probability $\hat{\mathcal{I}}(\omega_1)
\cup\mathcal{U}$ has a trifurcation point in $K$. For this, we first
observe from (\ref{eUgood}), (\ref{eqnuQ}) and (\ref{eqQK}) that
for some triple $\{z_1, z_2, z_3\} \subset K$,
\[
P_{z_1}^K \otimes P_{z_2}^K \otimes
P_{z_3}^K \bigl[\mbox{$\operatorname{Range}
\bigl(X^i\bigr)$, $i = 1,2,3$ are not connected through $\mathcal{U}$}
\bigr] > 0,
\]
where $X^i$, $i = 1,2,3$ stand for the random walk trajectories, which
are implicit in the above product space. The above can be rewritten as
\begin{eqnarray*}
&&
P_{z_1}^K \otimes P_{z_2}^K \otimes
P_{z_3}^K \bigl[ \mbox{$\operatorname{Range}
\bigl(X^1\bigr)$ and $\operatorname{Range}\bigl(X^2\bigr)$
are not connected}\\
&&\hspace*{80pt} \mbox{to each other
through $\mathcal{U}$ and $\operatorname{Range}
\bigl(X^3\bigr)$ does }\\
&&\hspace*{80pt}\mbox{not touch the clusters containing
$z_1$ and $z_2$ in}\\
&&\hspace*{152pt} \mbox{$\mathcal{U} \cup\operatorname {Range}
\bigl(X^1\bigr) \cup\operatorname{Range}\bigl(X^2\bigr)$}
\bigr] > 0.
\end{eqnarray*}
By conditioning in the range of the walks $X^1$ and $X^2$, we obtain that
\begin{eqnarray*}
&&
P_{z_1}^K \otimes P_{z_2}^K \otimes
P_{z_3}^K \otimes P_{z_3}^K \bigl[
\mbox{$\operatorname{Range}\bigl(X^i\bigr)$, $i = 1,2$ are not
connected through}\\
&&\hspace*{111pt}\mbox{$\mathcal{U}$ and $\operatorname{Range}\bigl(X^3
\bigr) \cup\operatorname{Range}\bigl(X^4\bigr)$ do not touch}\\
&&\hspace*{111pt}\mbox{the
clusters containing $z_1$ and $z_2$}\\
&&\hspace*{175pt}\mbox{in $\mathcal{U} \cup
\operatorname {Range}\bigl(X^1\bigr) \cup\operatorname{Range}
\bigl(X^2\bigr)$} \bigr]>0,
\end{eqnarray*}
and repeating this argument twice, we get
%
\begin{eqnarray}
\label{egianteq}
&&\bigl(P_{z_1}^K\bigr)^{\otimes2}
\otimes\bigl(P_{z_2}^K\bigr)^{\otimes2} \otimes
\bigl(P_{z_3}^K\bigr)^{\otimes2} \bigl[ \mbox{$\bigl(
\operatorname{Range}\bigl(X^1\bigr)\cup\operatorname{Range}
\bigl(X^2\bigr)\bigr)$,}\nonumber\\
&&\hspace*{139pt}\mbox{$\bigl(\operatorname{Range}
\bigl(X^3\bigr)\cup\operatorname{Range}\bigl(X^4\bigr)
\bigr)$,}\nonumber\\[-8pt]\\[-8pt]
&&\hspace*{139pt}\mbox{$\bigl(\operatorname{Range}\bigl(X^5\bigr)\cup
\operatorname{Range}\bigl(X^6\bigr)\bigr)$}\nonumber\\
&&\hspace*{139pt}\mbox{are not connected through
$\mathcal{U}$} \bigr]>0.\nonumber
\end{eqnarray}

We can construct a deterministic subgraph $G'$ of $K$ such that:
\begin{itemize}
\item$G'$ is connected,
\item$z_1$, $z_2$ and $z_3$ belong to $G'$ and
\item the removal of some vertex $y$ of $G'$ separates the vertices
$z_1$, $z_2$ and $z_3$ into disjoint components.
\end{itemize}
Indeed, to perform such construction, one could, for instance, take a
shortest path $\gamma$ from $z_1$ to $z_2$ in $K$ and a shortest path
$\gamma'$ from $z_3$ to the range of $\gamma$. Then, the graph $G'$
given by the union of the vertices and edges of $\gamma$ and $\gamma'$
would satisfy the above conditions with $y$ being the intersection of
$\gamma$ and $\gamma'$. Note that the three above conditions imply
that, for a subgraph $\mathcal{I}$ of $G$,
%
\begin{equation}
\label{etrifurc}
\begin{tabular}{p{323pt}}
if $\mathcal{I} \cap K = G'$, and
$z_1$, $z_2$ and $z_3$ belong to distinct
infinite connected
components of $\mathcal{I} \setminus(K \setminus\{
z_1,z_2,z_3\} )$, then $y$ is a
trifurcation point of $\mathcal{I}$.
\end{tabular}\hspace*{-27pt}
\end{equation}

Using\vspace*{1pt} that $G'$ is finite and connected, with positive
$P_{z_1}$-probability, the random walk $X^1$ can cover $G'$ and return
to $z_1$ before escaping from $K$. This together with (\ref{etrifurc})
and (\ref{egianteq}) gives that if $\mathcal{U}$ is good,
\[
\prod_{i = 1,2,3} \bigl(P_{z_i} \otimes
P_{z_i}^K\bigr) \biggl[\mbox{$y$ is a trifurcation point of }
\bigcup_{i = 1,\ldots, 6} \operatorname{Range}
\bigl(X^i\bigr) \cup\mathcal{U} \biggr]>0.
\]

Finally, we use (\ref{eqnuQ}) and (\ref{eqQK}) to conclude that if
$\mathcal{U}$ is good,
%
\begin{equation}
\mathbb{P}_1\bigl[\mbox{$y$ is a trifurcation point of $\hat{
\mathcal {I}}(\omega_1) \cup\mathcal{U}$}\bigr] > 0.
\end{equation}
This, together with Fubini's theorem and the fact that $\mathcal
{I}(\omega_2)$ is good with positive $\mathbb{P}_2$-probability, gives
that $y$ is a trifurcation point with positive $\mathbb
{P}_u$-probability. The rest of the proof of the theorem follows the
ideas of~\cite{BK89}; see also the proof of Theorem 2.4 in~\cite{HJ06}.
\end{pf}

\section{Domination by the frog model}\label{ssleepart}

Using the Poissonian character of random interlacements, we are going
to show that the cluster of $\mathcal{I}(\omega)$ containing a given
point $x \in V$ is dominated by the trace left by the particles in a
certain particle system, the so-called frog model. In what follows we
give an intuitive description of this process.

Place a random number of particles in each site of $V$ with a
prescribed distribution which will be given later. These particles
should be understood as ``sleeping'' when the system starts. In the
first iteration, only the particles sitting at the fixed site $x$ are
active, and they perform simple random walks. When an active particle
reaches the neighborhood of a sleeping one, the latter becomes active,
starting to perform simple random walk and so on. Note the two
differences between this model and the one studied in~\cite{RS04}: we
consider a random initial configuration instead of a deterministic one
and the active particles can wake a sleeping particle by simply
visiting its neighborhood, without the need to occupy the same site.\vadjust{\goodbreak}

We now give a precise description of this process, while keeping a
similar notation as that of previous sections. Consider the following
measure on the space $(W,\mathcal{W})$ of doubly infinite trajectories:
%
\begin{equation}
\label{eqmu} \tilde{\nu} = \sum_{y \in V}
{d_y} \tilde{P}_y,
\end{equation}
where $\tilde{P}_y$ is supported in $\{X_0 = y\} \subset W$ and given by
%
\begin{eqnarray}
\label{eqdeftilde}
&&
\tilde{P}_y \bigl[(X_{-i})_{i \geq0}
\in A, X_0 = y, (X_{i})_{i
\geq
0} \in B \bigr]
\nonumber\\[-8pt]\\[-8pt]
&&\qquad = P_y[A]P_y[B]\qquad \mbox{for every $A,B \in
\mathcal{W}_+$}.\nonumber
\end{eqnarray}
Intuitively speaking, the measure $\tilde{P}_y$ launches two
independent random walks from~$y$, one to the future and the other to
the past.

Consider also the following space of point measures on $W$:
%
\begin{eqnarray}
\label{eqM} M &=& \biggl\{\tilde{\omega} = \sum
_{i\geq 1} \delta_{w_i}; w_i \in W \mbox{ and
} \tilde{\omega} \bigl(\bigl\{X_0(w)\in K\bigr\} \bigr) < \infty,
\nonumber\\[-8pt]\\[-8pt]
&&\hspace*{135pt}\mbox{for all finite } K \subset V \biggr\},\nonumber
\end{eqnarray}
which is endowed with the sigma-algebra $\mathcal{M}$ generated by the
evaluation maps $\tilde{\omega} \to\tilde{\omega}(D)$, for any $D
\in
\mathcal{W}$.

In analogy to (\ref{eqinterlace}), for any $\tilde{\omega} \in M$
we define
%
\begin{equation}
\label{eqtildeinter} \tilde{\mathcal{I}} (\tilde{\omega}) = \biggl\{ \bigcup
_{w \in
\operatorname{supp}(\tilde{\omega})} \operatorname{Range}(w) \biggr\},
\end{equation}
which is the trace of the trajectories composing $\tilde{\omega}$.

We also introduce, in the probability space $(M,\mathcal{M},\tilde
{\mathbb{P}}_u)$, a Poisson point process $\tilde{\omega}$ in $W$ with
intensity measure given by\vspace*{1pt} ${u} \tilde{\nu}$. In Remark \ref
{rcoverall} below, we prove that, if $G$ has bounded degree, then
$\tilde{\mathcal{I}}(\tilde{\omega}) = V$, $\tilde{\mathbb{P}}_u$-a.s.
In our pictorial description, this corresponds to the fact that if we
wake up, all the particles at time zero all the sites will eventually
be visited. In what follows we will instead consider $\tilde{\mathcal
{I}} (\tilde{w}|_{\cdot})$ where~$\tilde{w}|_\cdot$ denotes the
restriction of $\tilde{w}$ to some set $\cdot\in\mathcal{W}$.

Fix now $\tilde{\omega} \in M$. The following construction can be
intuitively described as gradually ``revealing'' $\tilde{\omega}$. By
this, we mean that in each step we observe $\tilde{\omega}$ restricted
to larger and larger subsets of $W$. More precisely, recalling the
definition of $\overline K$ in Section~\ref{snotation}, define:
\begin{itemize}
\item$\tilde{A}_0(\tilde{\omega}) = \{x\}$ ``particles sitting at $x$
are activated at iteration $0$,'' and supposing we have constructed
$\tilde{A}_0, \ldots, \tilde{A}_{k-1}$,
\item let $\tilde{A}_k(\tilde{\omega}) = \overline{\tilde{\mathcal
{I}}(\tilde{\omega}  |_{\{X_0 \in\tilde{A}_{k-1}(\tilde{\omega
})\}
})}$, ``particles sitting at $\tilde{A}_{k-1}$ were activated and
performed random walks, whose ranges will determine the next active set
$\tilde{A}_k$.''\vadjust{\goodbreak}
\end{itemize}
Note that the above restriction of $\tilde{\omega}$ does not include
all trajectories which hit $\tilde{A}_{k-1}$, solely the ones starting
on it. Also, we sometimes write $\tilde{A}_k$ instead of $\tilde
{A}_k(\tilde{\omega})$ in order to avoid an overly heavy notation.

Due to the Poissonian character of $\tilde{\mathbb{P}}_u$, since
$\tilde
{A}_k(\tilde{\omega})$ is determined by $\tilde{\omega}|_{\{X_0 \in
\tilde{A}_{k-1}\}}$,
%
\begin{equation}
\label{eqinductilde}
\begin{tabular}{p{323pt}}
conditioned on $\tilde{A}_0(\tilde{
\omega}), \ldots, \tilde {A}_k(\tilde{\omega})$, the point measure $
\tilde{\omega}|_{\{X_0
\in
\tilde{A}_k \setminus\tilde{A}_{k-1}\}}$,
is distributed as a Poisson point process on $\{X_0 \in\tilde
{A}_k \setminus\tilde{A}_{k-1}\}$,
with intensity given by $u 1_{\{X_0 \in\tilde{A}_k \setminus
\tilde{A}_{k-1}\}} \cdot\tilde{\nu}$;
\end{tabular}\hspace*{-27pt}
\end{equation}
see, for instance,~\cite{R08}, Proposition 3.6.

We can now state our domination result.
%
\begin{proposition}
\label{pdomin}
Let $\mathcal{C}_x$ be the connected component of $\mathcal{I}(\omega)$
containing $x\in V$. We can find a coupling $Q_1$ between $\mathbb
{P}_u$ and $\tilde{\mathbb{P}}_u$ such that
%
\begin{equation}
\label{eqcouple1} \mathcal{C}_x \subset\tilde{\mathcal{I}} (
\tilde{\omega} |_{\{
X_0 \in\bigcup_k \tilde{A}_k\}} ),\qquad Q_1\mbox{-a.s.}
\end{equation}
\end{proposition}
\begin{pf} Recall the definition of $W^*_K$ above (\ref{eqomega}),
which we now extend to sets $K$ which are possibly infinite.

For a given $\omega\in\Omega$, we are going to construct sets
$A_0(\omega), A_1(\omega), \ldots$ in a similar way as we did for
$\tilde
{A}_k(\tilde{\omega})$'s. For this:
\begin{itemize}
\item consider $A_0 = \{x\}$ ``trajectories passing through $x$ are
activated at step $0$,'' and suppose we have constructed $A_0, \ldots, A_{k-1}$,
\item let $A_k(\omega) = \overline{\mathcal{I}(\omega
|_{W^*_{A_{k-1}}})}$ ``trajectories meeting $A_{k-1}$ were activated
and their ranges will determine the next active set $A_k$.''
\end{itemize}

Note that the Poissonian character of $\mathbb{P}_u$ gives us that for
every $k \geq0$,
%
\begin{equation}
\label{eqinducomega}
\begin{tabular}{p{323pt}}
conditioned on $A_0(\omega), \ldots,
A_k(\omega)$, the process $\omega|_{W^*_{A_k} \setminus W^*_{A_{k-1}}}$
is distributed as a Poisson point process on $W_{A_k}^*$ with
intensity
given by $u 1_{W^*_{A_k} \setminus W^*_{A_{k-1}}} \cdot\nu$;
\end{tabular}\hspace*{-27pt}
\end{equation}
compare with (\ref{eqinductilde}).

Another important remark is that, although $\mathcal{I}(\omega
|_{W^*_{\bigcup_k A_k}})$ could in principle be a proper subset of
$\mathcal{I}(\omega)$, actually the connected component $\mathcal{C}_x$
of $\mathcal{I}(\omega)$ containing $x$ is given by
%
\begin{equation}
\label{eqCxIsum} \mathcal{C}_x = \mathcal{I} ( \omega
|_{W^*_{\bigcup_k A_k}} ).
\end{equation}
To see\vspace*{1pt} why this holds, note first that $\mathcal{C}_x$ is empty if and
only if $\omega|_{W^*_{A_k}} = 0$ for every $k \geq0$. Assume now that
$\mathcal{C}_x$ is nonempty and take a point $y$ in $\mathcal{C}_x$,
which implies the existence of a path $x = x_0, x_1,\ldots, x_n = y$
contained in $\mathcal{I}(\omega)$. Supposing by contradiction that $y
\notin\mathcal{I} ( \omega|_{W^*_{\bigcup_k A_k}})$, let $j_o$ be the
first $j \leq n$ such that $x_j \notin\mathcal{I} ( \omega
|_{W^*_{\bigcup_k A_k}} )$. Since $\mathcal{C}_x$ is nonempty, we
conclude that $j_o \geq1$, and therefore $x_{j_o}$ is within distance
one from $x_{j_o-1} \in\mathcal{I}(\omega|_{W^*_{\bigcup_k A_k}})$. To
obtain a contradiction, let $k_o$ be such that $x_{j_o-1} \in\mathcal
{I}(\omega|_{W^*_{\bigcup_{k \leq k_o} A_k}})$ and observe that $x_{j_o}$
must belong to $A_{k_o + 1}$, and since it belongs to $\mathcal
{I}(\omega)$ it must also be in $\mathcal{I} ( \omega|_{W^*_{\bigcup_k
A_k}})$, which is a contradiction. This proves that $\mathcal{C}_x
\subseteq\mathcal{I} ( \omega|_{W^*_{\bigcup_k A_k}})$. The
other\vspace*{1pt}
inclusion in obvious since $\mathcal{I} ( \omega|_{W^*_{\bigcup_k A_k}})$
is connected and contained in $\mathcal{I}(\omega)$. This establishes
(\ref{eqCxIsum}).

From (\ref{eqCxIsum}), we see that the set $\mathcal{C}_x$ can be
written in a way that resembles the set $\tilde{\mathcal{I}}( \tilde
{\omega}|_{\{X_0 \in\bigcup_k \tilde{A}_k\}})$ appearing in (\ref
{eqcouple1}). But to proceed with the proof of the proposition, we
first need to find a coupling $Q_1$ between $\omega$ (under $\mathbb
{P}_u$) and $\tilde{\omega}$ (under $\tilde{\mathbb{P}}_u$) such that
%
\begin{equation}
\label{eqdomomegas} \omega |_{W^*_{A_k(\omega)}} \leq\pi^* \circ (\tilde {\omega}
|_{\{X_0 \in A_k(\omega)\}} )\qquad \mbox{for every $k \geq 0$, $Q_1$-a.s.}
\end{equation}
Note that the sets $\tilde{A}_k(\tilde{\omega})$ do not appear in the
above equation. In order for $Q_1$ to satisfy the above, it is clearly
enough that
%
\begin{equation}
\label{edominter} \qquad\omega |_{W^*_{A_k} \setminus W^*_{A_{k-1}}} \leq\pi^* \circ ( \tilde{\omega}
|_{\{X_0 \in A_k \setminus A_{k-1}\}} )\qquad\mbox{for every $k \geq0$,
$Q_1$-a.s.},
\end{equation}
where we used the convention that $A_{-1} = \varnothing$. We now prove
the existence $Q_1$ satisfying (\ref{edominter}) and consequently
(\ref{eqdomomegas}).

Recall first that, conditioned on $A_1, A_2, \ldots, A_k$, the left-hand
side of inequality (\ref{edominter}) is independent of $\omega
|_{W^*_{A_{k-1}}}$ while the right-hand side is independent of $\tilde
{\omega}|_{\{X_0 \in A_{k-1}\}}$. Therefore, to be able to construct
the coupling $Q_1$ satisfying (\ref{edominter}), it suffices to show
that for any $B \subseteq B'$,
%
\begin{equation}
\label{edomBB}
\begin{tabular}{p{323pt}}
we can couple $\omega |_{W^*_{B'} \setminus W^*_{B}}$ with $\tilde{
\omega} |_{\{X_0 \in B' \setminus B\}}$ in a way that\\
\mbox{$\omega |_{W^*_{B'} \setminus W^*_{B}} \leq\pi^* \circ ( \tilde {\omega} |_{\{X_0 \in B' \setminus B\}} )$}.
\end{tabular}\hspace*{-27pt}
\end{equation}
In fact, once we have established the above, we can proceed by
induction with $B' = A_k$ and $B = A_{k-1}$, for $k \geq0$ and define
$\omega |_{W^* \setminus\bigcup_k W^*_{A_k}}$ and $\tilde{\omega
}
|_{\{X_0 \notin\bigcup_k A_k\}}$ independently. This way, they will have
the correct marginal distribution (see, e.g., Proposition 3.6
in~\cite{R08}) and will satisfy (\ref{edominter}).

As a further reduction, we claim that it is enough to establish (\ref
{edomBB}) in the case where $B' = B \cup\{y\}$, with $y \notin B$.
Indeed, if $B' = B \cup\{y_1, y_2, \ldots\}$, then we write $B_i = B
\cup\{y_1, \ldots, y_i\}$ and use (\ref{edomBB}) repeatedly for the
sets $B_{i+1}$ and $B_i$, with $i \geq0$ to obtain (\ref{edomBB}).

From now on, fix $B \subset V$ and $y \notin B$ and note, by (\ref
{eqnuQ}), that
%
\begin{equation}
\label{eintBy} \mbox{the intensity measure of $\omega |_{W^*_{B \cup\{y\}}
\setminus W^*_{B}}$ is
given by $u \pi^* \!\circ\! (\mathbf{1}_{{\{H_B
= \infty\}}} \!\cdot\! Q_{\{y\}} )$,}\hspace*{-35pt}
\end{equation}
which we estimate as follows.

Fix $C, C' \in\mathcal{W}_+$ and consider the event $D \in\mathcal
{W}$ given by $D = \{H_y = 0\} \cap \{(X_{-i})_{i \geq0} \in C
\} \cap \{ (X_{i})_{i \geq0} \in C'  \}$. Then
%
\begin{eqnarray}
&&
\mathbf{1}_{{\{H_B = \infty\}}} \cdot Q_{\{y\}} (D) = Q_{\{
y\}} [D,
H_B = \infty]
\nonumber
\\
&&\qquad \stackrel{(\rref{eqQK})} {=} P_{y} [C, H_{B} =
\tilde{H}_y = \infty]{d_y } P_{y}
\bigl[C', H_{B} = \infty\bigr]
\nonumber
\\
&&\qquad \le P_{y} [C, \tilde{H}_y = \infty]{d_y
} P_{y} \bigl[C'\bigr]
\\
&&\qquad \stackrel{(\rref{eqdeftilde})} {=} {d_y\,}\tilde{P}_y
\bigl[ (X_{-i})_{i \geq0} \in C \cap\{\tilde{H}_y =
\infty\}, X_0 = y, (X_i)_{i \geq0} \in
C' \bigr]
\nonumber
\\
&&\qquad = {d_y\,}\tilde{P}_y[D].
\nonumber
\end{eqnarray}
Since the above holds for every event $D$ as above, we conclude that
$\mathbf{1}_{{\{H_B = \infty\}}} \cdot Q_{\{y\}} \leq d_y \tilde{P}_y$.
We can therefore construct a Poisson point process $\tilde{\omega}_-$
in $W$ with intensity given by $u \mathbf{1}_{{\{H_B = \infty\}}}
\cdot
Q_{\{y\}}$ in a way that $\tilde{\omega}_- \leq\tilde{\omega}_{\{
X_0 =
y\}}$. Since $\pi^* \circ\tilde{\omega}_-$ has the same law as
$\omega
|_{W^*_{B \cup\{y\}} \setminus W^*_{B}}$, we conclude (\ref
{edomBB}) from (\ref{eintBy}), for the case where $B' = B \cup\{y\}$.
As we have discussed, this is enough to establish (\ref{edomBB}) in
general and consequently (\ref{edominter}) and (\ref{eqdomomegas}).

Now that we have constructed $Q_1$, we are going to prove that it
satisfies (\ref{eqcouple1}). First we claim that
%
\begin{equation}
\label{eqAtildeA} A_k(\omega) \subseteq\tilde{A}_k(
\tilde{\omega})\qquad \mbox{for every $k \geq0$, $Q_1$-a.s.}
\end{equation}
To prove that, observe first that $A_0 = \tilde{A}_0 = \{x\}$ and
suppose that we have established the above result for $k-1$. Then we
use (\ref{eqdomomegas}) to obtain that, $Q_1$-almost surely,
%
\begin{equation}
\mathcal{I} ( \omega |_{W^*_{A_{k-1}}} ) \subseteq \tilde {\mathcal{I}} (\tilde{
\omega} |_{\{X_0 \in A_{k-1}\}} )
\end{equation}
and therefore
%
\begin{equation}
A_k(\omega) ={} \overline{\mathcal{I} ( \omega |_{W^*_{A_{k-1}}} )}
\subseteq\overline{\tilde{\mathcal{I}} (\tilde{\omega} |_{\{X_0 \in A_{k-1}\}} )}
\subseteq \overline{\tilde{\mathcal{I}} (\tilde{\omega} |_{\{X_0
\in
\tilde{A}_{k-1}\}} )} =
\tilde{A}_k,
\end{equation}
proving (\ref{eqAtildeA}).

Finally, using (\ref{eqCxIsum}), (\ref{eqdomomegas}) and (\ref
{eqAtildeA}), we obtain
%
\begin{equation}\qquad
\mathcal{C}_x = \mathcal{I}(\omega|_{W^*_{\bigcup_k A_k}}) \subseteq
\tilde {\mathcal{I}}( \tilde{\omega}|_{\{X_0 \in\bigcup_k A_k\}}) \subseteq \tilde{
\mathcal{I}}(\tilde{\omega}|_{\{X_0 \in\bigcup_k \tilde
{A}_k\}}) ,\qquad\mbox{$Q_1$-a.s.},
\end{equation}
proving (\ref{eqcouple1}).
\end{pf}
%
\begin{remark}
\label{rcoverall}
As promised below (\ref{eqtildeinter}), let us show that if $G$ has
degree bounded by $\Delta$, then $\tilde{\mathbb{P}}_u$-a.s. $\tilde
{\mathcal{I}}(\tilde{\omega}) = V$. First, fix some $y \in V$ and
consider the events $C_z = \{$some walker started at $z$ hits $y\}$,
for $z \in V$. Clearly, under the measure $\tilde{\mathbb{P}}_u$ they
are independent. We now estimate $\sum_{z\in V} \tilde{\mathbb{P}}_u
[C_z]$, which is bounded from below by
%
\begin{eqnarray}
&&
\sum_{z\in V}  \tilde{\mathbb{P}}_u
\bigl[\tilde{\omega}(X_0 \in z) > 0\bigr] P_z[H_y
< \infty]\nonumber\\
&&\qquad \geq c_u \sum_{z\in V}
P_z[H_y < \infty]
\nonumber\\[-8pt]\\[-8pt]
&&\qquad = c_u \frac{1}{E_y[\mathrm{number}\ \mathrm{of}\ \mathrm{visists}\ \mathrm{to}\ y]} \sum_{z\in V}
E_z[\mathrm{number}\ \mathrm{of}\ \mathrm{visists}\ \mathrm{to}\ y]
\nonumber\\
&&\qquad \geq c_{u,y} \sum_{z \in V} \sum
_{n \geq0} p^{(n)}(z,y) \geq c_{u,y,
\Delta} \sum
_{n \geq0} \sum_{z \in V}
p^{(n)}(y,z) = \infty,
\nonumber
\end{eqnarray}
where in the second line we used the strong Markov property, and in the
third line we used the reversibility of the walk.

Using the independence of the events $\{C_z\}_{z \in V}$ and the
Borel--Cantelli lemma, we conclude that $\tilde{\mathbb{P}}_u$-a.s. the
point $y$ is visited by infinitely many walks. Since this holds true
for any of the countable $y$'s in $V$, we conclude that $\tilde
{\mathbb
{P}}_u$-a.s. $\tilde{\mathcal{I}}(\tilde{\omega}) = V$, as desired.
\end{remark}

\section{Domination by branching random walk}\label{sbranching}

In this section we show that the so-called frog model defined in
Section~\ref{ssleepart} can be dominated by a certain multitype
branching random walk which we introduce below. Although the main
result of this section is very intuitive, its proof is somewhat
technical. This is mainly due to the infection of neighbors present in
our version of the frog model (see comment in the second paragraph of
Section~\ref{ssleepart}) and the notation we chose for the frog model,
resembling random interlacements' formalism.

Consider two types of individuals, namely $\circ$ and $\bullet$, and
four probability distributions on $\mathbb{Z}_+$: $q_{\circ
\shortrightarrow\circ}, q_{\circ\shortrightarrow\bullet},
q_{\bullet
\shortrightarrow\circ}$ and $q_{\bullet\shortrightarrow\bullet}$.
These are the offspring distributions which should be understood as
follows: $q_{\circ\shortrightarrow\bullet}$ is the distribution of
the number of $\bullet$-descendants of a $\circ$-individual, and
analogously for $q_{\circ\shortrightarrow\circ}, q_{\bullet
\shortrightarrow\circ}$ and $q_{\bullet\shortrightarrow\bullet}$. If
one supposes that all individuals have an independent number of $\circ$
and $\bullet$ descendants, then these four distributions completely
characterize the branching mechanism.

We now make a specific choice for the $q$'s, namely:
\begin{itemize}
\item$q_{\circ\shortrightarrow\bullet} = \delta_2$, ``$\circ
$-individuals always have two descendants of type $\bullet$,''
\item$q_{\bullet\shortrightarrow\bullet} = \delta_1$, ``$\bullet
$-individuals always have one descendants of type $\bullet$'' and
\item$q_{\circ\shortrightarrow\circ} = q_{\bullet\shortrightarrow
\circ} = \operatorname{Poisson}({u \Delta^2})$.
\end{itemize}
Finally we need to introduce the initial distribution:
%
\begin{equation}
\begin{tabular}{p{323pt}}
at the first generation, all individuals are of type $\circ$
and their number is $\operatorname{Poisson}({u \Delta^2}
)$ distributed.
\end{tabular}\hspace*{-27pt}
\end{equation}

We now construct in some abstract probability space $(S, \mathcal{S},
\bfPu)$ the tree $\mathcal{T} = (T,\mathcal{E})$ originated by the
above $2$-type Galton--Watson process, which is completely
characterized by the above properties. Rigorously, $\mathcal{T}$ should
be called a forest, since it has as many connected components as the
number of individuals at generation one. However, we keep the notation
``tree'' for simplicity.

\begin{figure}

\includegraphics{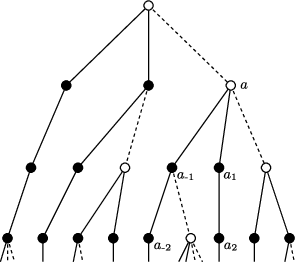}

\caption{The multitype Galton--Watson tree $\mathcal{T} = (T,\mathcal
{E})$, where $T = T^\circ\cup T^\bullet$.}
\label{ftree}
\end{figure}

Here, the set of vertices $T$ is given by the disjoint union of
$T^\circ
$ and $T^\bullet$, corresponding to the type of each individual; see
Figure~\ref{ftree}. Supposing that the degree of $G$ is bounded by
$\Delta$, we construct in the same probability space $(S, \mathcal{S},
\bfPu)$ a collection $(\ell_e)_{e \in\mathcal{E}}$ of i.i.d. random
variables with uniform distribution over the set $\{1, 2, \ldots,
\Delta
!\}$. These variables $\ell_e$ should be regarded as random labels
assigned to the edges of $\mathcal{T}$.

We now use the above construction to define a special branching random
walk on the original graph $G = (V,E)$. For this, for every vertex $y
\in V$, we fix a bijection $\phi_y\dvtx (\mathbb{Z}/d_y\mathbb{Z}) \to
\mathcal{N}(y)$. Recall that $d_y$ is the degree of $y$, and $\mathcal
{N}(y)$ stands for the set of its neighbors.

Fix now $x \in V$, a tree $\mathcal{T}$, the labels $(\ell_e)_{e \in
\mathcal{E}}$ and $(\phi_y)_{y \in V}$ as above. We define the
coordinate processes $(Z_a^x)_{a \in T}$ of the branching random walk by:
%
\begin{equation}
\begin{tabular}{p{323pt}}
$Z_a^x$ = $x$ for every $a$ in the first
generation of $\mathcal {T}$.
If $Z_a^x = y$ for some $a \in T$, then for any
descendant
$b$ of $a$, we define $Z_b^x$ as $
\phi_y(\ell_{\{a,b\}} \operatorname{mod} d_y)$.
\end{tabular}\hspace*{-27pt}
\end{equation}

Recall that $\ell_e$ was uniformly chosen in $\{1, 2, \ldots, \Delta
!\}
$, and note that $d_y$ divides~$\Delta!$. This implies that,
conditioned on $Z_a^x = y$, $Z_b^x$ is a uniformly chosen neighbor of
$y$, for every $b$ descendant of $a$ in $\mathcal{T}$. Using this fact
inductively, we conclude that given $\mathcal{T}$, for any $b \in T$,
%
\begin{eqnarray}
\begin{tabular}{p{323pt}}
\label{ewalk}
$Z_b^x$ under the probability
distribution $\bfPu$ has the same law as the
random walk $X_n$ under $P_x$, where $n$ is the
generation of $b$ in $\mathcal{T}$.
\end{tabular}\hspace*{-27pt}
\end{eqnarray}

From our specific choice of the offspring distribution, note that
$\bfPu$-a.s.
%
\begin{equation}
\begin{tabular}{p{323pt}}
for every point $a \in T^\circ$, there exist exactly two
infinite
lines of \mbox{$\bullet$-descendants} going down from $a$.
\end{tabular}\hspace*{-27pt}
\end{equation}
We denote the union of these two lines by $L_a \subset{\mathcal T}$,
which are drawn with continuous segments in Figure~\ref{ftree}. Denote
by $E(L_a)\subset{\mathcal E}$ the set of edges in $L_a$ and by
$V(L_a)\subset T$ the set of individuals in $L_a$.

Having constructed the branching random walk $(Z_a^x)_{a \in T}$, we
can use the above observation to associate, for every $a \in T^\circ$ a
doubly infinite trajectory $w_a \in W$ given by the image of $L_a$
under the branching random walk $Z$. And using (\ref{ewalk}), we
conclude that
%
\begin{equation}
\label{eqindtraj} \mbox{the law of $w_a$ as an element of $W$ is
given by $\tilde {P}_{Z_{a}}$ as in (\ref{eqdeftilde}).}
\end{equation}

The above allows us to define
\[
\mathring{\omega} = \sum_{a \in T^\circ}
\delta_{w_a} \in M.
\]

We can now prove the domination result, which gives a way to control
the particle system defined in the previous section with the above
branching random walk.

%
\begin{proposition}
\label{pdomin2} If the degree of $G$ is bounded by $\Delta$, there is
a coupling $Q_2$ between the laws $\tilde{\mathbb P}_u$ and $\bfPu$,
such that $Q_2$-a.s. we have
%
\begin{equation}
\tilde{\mathcal{I}} ( \tilde{\omega} |_{\{X_0 \in\bigcup_k
\tilde
{A}_k\}} ) \subseteq\tilde{
\mathcal{I}} ( \mathring{\omega } ) = \bigl\{Z_a^x; a
\in T\bigr\}.
\end{equation}
\end{proposition}
\begin{pf}
We recall the definition of the sets $\tilde{A}_k(\tilde\omega)$
constructed above Proposition~\ref{pdomin}. We will need an analogous
way to gradually reveal the measure $\mathring{\omega}$ which can be
informally described as ``slowly revealing the lines $L_a$'' composing
$\mathcal{T}$. More precisely:
\begin{itemize}
\item consider $\mathtt{A}_0 = T_0^\circ$, ``all individuals in the
first generation of $\mathcal{T}$'' and supposing we have constructed
$\mathtt{A}_{k-1}$,
\item let $L_{\mathtt{A}_{k-1}}$ be the union of all $L_a$ where $a
\in
\mathtt{A}_{k-1}$ and all edges in $\mathcal{T}$ that connect any such
$a$ to its parent. Then define $\mathtt{A}_k$ to be the union of
$V(L_{\mathtt{A}_{k-1}})$ with its descendants.
\end{itemize}
In Figure~\ref{fadd}, we show the sets $L_{\mathtt{A}_0}$,
$L_{\mathtt
{A}_1}$ and $L_{\mathtt{A}_2}$.

\begin{figure}

\includegraphics{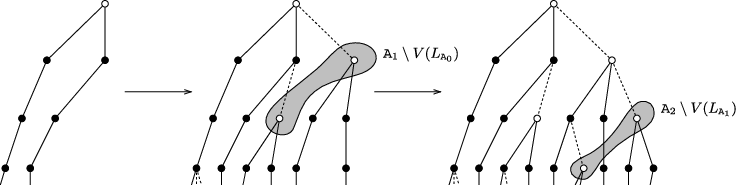}

\caption{The sets $L_{\mathtt{A}_0}$, $L_{\mathtt{A}_1}$ and
$L_{\mathtt{A}_2}$ corresponding to $\mathcal{T}$ in Figure \protect
\ref{ftree}.}
\label{fadd}
\end{figure}

Observe also the following two consequences of our particular choice
for the offspring distribution of $\mathcal{T}$:
\begin{itemize}
\item$\mathtt{A}_k \setminus V(L_{\mathtt{A}_{k-1}})$ consists solely
of $\circ$-type individuals,
\item conditioned on $L_{\mathtt{A}_{k-1}}$, the set $\mathtt{A}_k
\setminus V(L_{\mathtt{A}_{k-1}})$ can be obtained by independently
attaching a $\operatorname{Poisson}(u \Delta^2)$ number of $\circ$-offsprings
to each site of $L_{\mathtt{A}_{k-1}} \setminus L_{\mathtt{A}_{k-2}}$,
for any $k \geq1$,
\end{itemize}
where we use the convention that
$L_{\mathtt{A}_{-1}}=\varnothing$.\vadjust{\goodbreak}

In order to compare $\tilde{\omega}$ and $\mathring{\omega}$, we
claim that
%
\begin{equation}
\label{eqbpoisson}
\begin{tabular}{p{303pt}}
given $L_{\mathtt{A}_{k-1}}$ and $(\ell_e)_{e \in E(L_{\mathtt
{A}_{k-1}})}$,
let ${\mathcal Z}_k = (Z_a^x
)_{a\in V(L_{\mathtt
{A}_{k-1}})}$.
Then ${\sum_{a \in{\mathtt{A}_k}\setminus V(L_{\mathtt
{A}_{k-1}})} \delta_{w_a}}$ is
a Poisson point process
with intensity measure larger or equal than
${\sum
_{x \in
\overline{{\mathcal Z}}_k \setminus{\mathcal Z}_k} u d_x \tilde P_x}$.
\end{tabular}\hspace*{-27pt}
\end{equation}

To prove the above statement, observe first that ${\mathcal Z}_k$ can
be written in terms of $L_{\mathtt{A}_{k-1}}$ and $(\ell_e)_{e \in
E(L_{\mathtt{A}_{k-1}})}$. By using the comment above (\ref
{eqbpoisson}), we conclude that $(Z_a^x)_{a \in\mathtt{A}_k
\setminus
V(L_{\mathtt{A}_{k-1}})}$ is a Poisson process in $\overline
{{\mathcal
Z}}_k$ with intensity given by
%
\begin{equation}
\label{eqintAk} \sum_{a \in V(L_{\mathtt{A}_{k-1}}) \setminus V(L_{\mathtt{A}_{k-2}})} \frac{{u \Delta^2}}{\deg(Z_a^x)} \sum
_{y \in\mathcal{N}(Z_a^x)} \delta_y.
\end{equation}
One can now conclude (\ref{eqbpoisson}) from (\ref{eqintAk}) and
(\ref{eqindtraj}) and $u \Delta^2/\deg(Z_a^x)$ $\ge$ $u d_x$.

We now prove that
%
\begin{equation}
\label{eqdomin3}
\begin{tabular}{p{323pt}}
there exist a coupling $Q_2$ between $\tilde{
\mathbb{P}}_u$ and $\mathbf{P}_u$,
such that $\tilde{\omega}|_{\{X_0 \in\tilde{A}_k\}} \leq \sum
_{a \in\mathtt{A}_k} \delta_{w_a}$, $Q_2$-a.s.
for every $k \geq0$.
\end{tabular}\hspace*{-27pt}
\end{equation}

For this,\vspace*{1pt} fix a point measure $\tilde\omega$ sampled according to
$\tilde{\mathbb{P}}_u$. This gives rise to a sequence of sets $\tilde
{A}_0$, $\tilde{A}_1,\ldots\,$. Note that $\tilde{\omega}(\{X_0 \in
\tilde{A}_0\})$ is Poisson distributed, with a parameter not bigger
than $u \Delta$. We can therefore couple this random variable with the
initial number $|T_0|$ of individuals in the Galton--Watson tree, in a
way that almost surely
\[
\tilde{\omega}\bigl(\{X_0 \in\tilde{A}_0\}\bigr)
\leq|T_0|.
\]
And on the above event, one can construct the coupling in (\ref
{eqdomin3}) for $k = 0$, using~(\ref{eqindtraj}).

We now proceed by induction. Suppose that we have obtained the coupling
between $\tilde{\omega}|_{\{X_0 \in\tilde{A}_{k-1}\}}$ under
$\tilde
{\mathbb{P}}_u$ and $L_{\mathtt{A}_{k-1}}$ under $\mathbf{P}_u$
[together with its labels $(\ell_e)_{e \in E(L_{\mathtt{A}_{k-1}})}$]
in a way that $\tilde{\omega}|_{\{X_0 \in\tilde{A}_{k-1}\}}$ $\leq$
$\sum_{a \in\mathtt{A}_{k-1}} \delta_{w_a}$ for a given $k$. Then
$\tilde{A}_{k} \subseteq\overline{\mathcal{Z}}_k$. We can now obtain
the statement for $k$ using (\ref{eqbpoisson}) and comparing the
intensity measure of the referred Poisson point process with that of
(\ref{eqmu}).

The statement in Proposition~\ref{pdomin2} clearly follows from
(\ref{eqdomin3}).
\end{pf}

\section{The disconnected phase}\label{sdisconnsection}

The aim of this section is to show that on any nonamenable graph of
bounded degree, $u_c > 0$; that is, for some $u>0$ small enough the
interlacement set is $\mathbb{P}_u$-a.s. disconnected; see Theorem
\ref
{tdisconn}. First, in Section~\ref{sbrwres}, we obtain estimates
for the branching random walk introduced in Section~\ref{sbranching},
and then those estimates are used in Section~\ref{smthmproof} to
establish the existence of such a disconnectedness phase.

\subsection{Results for the branching random walk}\label{sbrwres}
The goal of this subsection is to prove a heat-kernel estimate for the
branching random walk; see Proposition~\ref{tbrwthm} below. In
Section~\ref{smthmproof} this will be used as a key ingredient for
showing Theorem~\ref{tdisconn}. Recall the definitions of ${\mathcal
T}=(T,\mathcal E)$, $T^\circ$ and $T^\bullet$ from Section \ref
{sbranching}. We introduce some additional notation.
For $n\in{\mathbb N}$, let
%
\begin{eqnarray}
\label{egendef}T_n&=&\{a\in T; a \mbox{ is in generation } n\},
\\
%
\label{egendef1}T_n^\circ&=&\bigl
\{a\in T^\circ; a \mbox{ is in generation } n\bigr\}
\end{eqnarray}
and
%
\begin{equation}
\label{egendef2}T_n^\bullet=\bigl\{a\in T^\bullet; a
\mbox{ is in generation } n\bigr\}.
\end{equation}
Then clearly $|T_n|=|T_n^\circ|+|T_n^\bullet|$.

We begin by bounding the expected number of members of $T_n$:
%
\begin{lemma}\label{lexpmemberbound} For every $n \ge0$ and $u \ge0$,
%
\begin{equation}
\label{eexpmembern} \E_u\bigl[|T_n|\bigr]\le c_{u,\lambda}
\bigl(1+2u \Delta^2\bigr)^{n}.
\end{equation}
\end{lemma}
\begin{pf}
To show (\ref{eexpmembern}), we proceed by fairly standard recursion
arguments. First, using the fact that all members in generation $0$ are
of type $\circ$, we observe that
%
\begin{equation}
\label{erec1} \E_{u}\bigl[|T_0|\bigr]=\E_{u}
\bigl[\bigl|T_0^\circ\bigr|\bigr] = u \Delta^2
\end{equation}
and
%
\begin{equation}
\label{erec11} \E_{u}\bigl[|T_1|\bigr]=\bigl(2+u
\Delta^2\bigr)\E_{u}\bigl[\bigl|T_0^\circ\bigr|
\bigr]=\bigl(2+u \Delta^2\bigr) u \Delta^2.
\end{equation}
In addition, we have
%
\begin{eqnarray}
\label{erec2} \E_u\bigl[\bigl|T_n^\circ\bigr|\bigr] &=& u
\Delta^2\E_u\bigl[\bigl|T_{n-1}^\circ\bigr|
\bigr]+u \Delta^2\E_u\bigl[\bigl|T_{n-1}^\bullet
\bigr|\bigr]
\nonumber\\[-8pt]\\[-8pt]
&=& u \Delta^2\E_u\bigl[|T_{n-1}|\bigr]
\qquad\mbox{for every }n\ge1.\nonumber
\end{eqnarray}
Hence, for every $n\ge2$,
%
\begin{eqnarray}
\label{erec3} \E_u\bigl[|T_n|\bigr]& = &\bigl(1+u
\Delta^2\bigr)\E_u\bigl[\bigl|T_{n-1}^\bullet\bigr|
\bigr]+\bigl(2+u \Delta^2\bigr)\E_u\bigl[\bigl|T_{n-1}^\circ
\bigr|\bigr]
\nonumber
\\
& = &\bigl(1+u \Delta^2\bigr)\E_u\bigl[|T_{n-1}|\bigr]+
\E_u\bigl[\bigl|T_{n-1}^\circ\bigr|\bigr]
\nonumber
\\
& \stackrel{(\rref{erec2})} {=}& \bigl(1+u \Delta^2\bigr)
\E_u\bigl[|T_{n-1}|\bigr]+u \Delta^2\E_u\bigl[|T_{n-2}|\bigr]
\\
& \le &\bigl(1+u \Delta^2\bigr)\E_u\bigl[|T_{n-1}|\bigr]+u
\Delta^2\E_u\bigl[|T_{n-1}|\bigr]
\nonumber
\\
& \stackrel{ \phantom{(\ref{erec2})} } {=} &\bigl(1+2u \Delta^2\bigr)
\E\bigl[|T_{n-1}|\bigr].\nonumber
\end{eqnarray}
Now (\ref{eexpmembern}) follows from (\ref{erec3}), (\ref{erec11})
and induction.
\end{pf}

The proposition below can be seen as a heat-kernel estimate for the
above constructed branching random walk.
%
\begin{proposition}\label{tbrwthm}
Suppose $G=(V,E)$ is nonamenable, with degree bound\-ed by $\Delta<
\infty$. Recall the definition of the branching random walk starting at
a $x\in V$, introduced in Section~\ref{sbranching}. There is $u_0>0$
such that, for every $x \in V$ fixed,
%
\begin{equation}
\label{ehitconv} \lim_{y; d(x,y)\to\infty}\P_{u}\bigl[\exists a\in T\dvtx
Z_a^x=y\bigr]=0\qquad\mbox{for every }u<u_0.
\end{equation}
\end{proposition}
\begin{pf}
By~\cite{W00}, Lemma (8.1), page 84, we have
%
\begin{equation}
\label{espektral} p^{(n)}(x,y)\le c \rho^n\qquad\mbox{for every
$x,y \in V$.}
\end{equation}
Denote the event appearing in the probability on the left-hand side
of (\ref{ehitconv}) with $H_y$ (where the index $x$ is omitted since
it is kept fixed during the proof).

We observe the following inclusion:
%
\begin{equation}
\label{eevincl} H_y \subset\bigcup_{n=d(x,y)}^{\infty}
\bigl\{\bigl|\bigl\{a\in T_n\dvtx  Z_a^x=y\bigr\} \bigr|
\ge 1\bigr\}
\end{equation}
and suppose that
%
\begin{equation}
\label{eusmall} u<\frac{1}{2\Delta^2} \biggl(\frac{1}{\rho}-1 \biggr).
\end{equation}
We now obtain
\begin{eqnarray*}
\P_{u}[H_y] & \stackrel{(\rref{eevincl})} {\le} &\sum
_{n=d(x,y)}^{\infty}\P_{u}\bigl[\bigl|\bigl\{a
\in T_n\dvtx  Z_a^x=y\bigr\} \bigr|\ge1\bigr]
\\
&\le&\sum_{n=d(x,y)}^{\infty}\E_{u}
\bigl[\bigl|\bigl\{a\in T_n\dvtx  Z_a^x=y\bigr\} \bigr|
\bigr]
\\
& = &\sum_{n=d(x,y)}^{\infty} \sum
_{m \geq1} \E_{u}\bigl[\bigl|\bigl\{a\in T_n\dvtx
Z_a^x=y\bigr\}\bigr|, |T_n| = m\bigr]
\\
& = &\sum_{n=d(x,y)}^{\infty} \sum
_{m \geq1} \E_{u} \biggl[\sum
_{a \in T_n} \E_{u}\bigl[Z_a^x=y|
\mathcal{T}\bigr],|T_n| = m \biggr].
\end{eqnarray*}
By using the independence of $\mathcal{T}$ and the labels
$(\ell_e)_{e \in\mathcal{E}}$, we deduce from (\ref{eexpmembern}) and
(\ref{ewalk}) that the above is bounded by
\begin{eqnarray*}
& \stackrel{} {\le} &\sum_{n=d(x,y)}^{\infty}c_{u, \lambda}
\bigl(1+2u \Delta^2\bigr)^{n} P_x[X_n=y]
\\
&\stackrel{(\rref{espektral})} {\le} &\sum_{n=d(x,y)}^{\infty
}c_{u, \lambda}
\bigl(\rho\bigl(1+2u \Delta^2\bigr) \bigr)^{n},
\end{eqnarray*}
which converges to zero as $d(x,y)$ goes to infinity, due to (\ref{eusmall}).

Therefore, we see that (\ref{ehitconv}) holds with
%
\begin{equation}
\label{euchoice} u_0=\frac{1}{2\Delta^2} \biggl(\frac{1}{\rho}-1
\biggr)>0,
\end{equation}
completing the proof of the theorem.
\end{pf}

\subsection{Disconnectedness when $u$ is small}\label{smthmproof}

We now collect the results of previous sections in order to obtain the
almost sure disconnectedness of random interlacements at low levels for
a large class of nonamenable graphs of bounded degree. Before, recall
the definitions of the spectral radius $\rho$ in (\ref{spectraldef})
and the capacity in~(\ref{eqcapeq}).
%
\begin{theorem}\label{tdisconn}
Consider a nonamenable graph $G$ with degree bounded by $\Delta<
\infty$. Then
%
\begin{equation}
\label{mainstat} u_c\ge\frac{1}{2 \Delta^2} \biggl(\frac{1}{\rho}-1
\biggr).
\end{equation}
\end{theorem}
\begin{pf}
We first observe that the nonamenability, together with the
assumption that the degrees are uniformly bounded imply that
%
\begin{equation}
\label{eloweresc} \inf_{x \in V} P_x[\tilde{H}_x=
\infty] >0\qquad \mbox{see, for instance, (4.4) of~\cite{Tei09}.}
\end{equation}

For $u,v\in V$, write $u\stackrel{\lacei(\omega
)}{\longleftrightarrow}
v$ if $u$ and $v$ belong to the same connected component of the
subgraph induced by ${\mathcal I}(\omega)$. If $u>0$ is such that
$\lacei(\omega)$ is connected $\hP_u$-a.s., then for $o \in V$,
%
\begin{eqnarray}
\label{conneq} \inf_{y\in V}\hP_u \bigl[o\stackrel{\lacei(
\omega )} {\longleftrightarrow} y \bigr] &=& \inf_{y\in
V}\hP_u
\bigl[o\in\lacei(\omega), y\in\lacei(\omega)\bigr]
\nonumber
\\
&\ge&\inf_{y\in
V}\hP_u\bigl[o\in\lacei(\omega)\bigr]
\hP_u\bigl[y\in\lacei(\omega)\bigr]
\nonumber
\\
&\ge &\inf_{y\in V}\bigl(1-\exp\bigl\{-u d_y
P_y[\tilde{H}_y=\infty]\bigr\}\bigr)^2
\\
&\ge&\Bigl(1-\exp\Bigl\{-u \inf_{x \in V} P_x[
\tilde{H}_x=\infty]\Bigr\}\Bigr)^2
\nonumber
\\
&>&0,\nonumber
\end{eqnarray}
where the FKG-inequality (see~\cite{Tei09}, Theorem 3.1) was
used in the first inequality, and in the third inequality we used the
fact that $d_x\ge1$. Consequently, to show that with positive
probability $\lacei(\omega)$ is disconnected
for small $u$, it is enough to show that
\[
\lim_{d(o,y)\to\infty}\hP_u \bigl[o\stackrel{\lacei(\omega )} {
\longleftrightarrow} y \bigr]=0,
\]
when $u$ is sufficiently small. From Propositions~\ref{pdomin}
and~\ref{pdomin2} we know that
%
\begin{equation}
\label{edominancecons} {\mathbb P}_u \bigl[o\stackrel{\lacei(
\omega)} {\longleftrightarrow} y \bigr]\le{\mathbf P}_u[H_y],
\end{equation}
where we recall the definition of the event $H_y$ from the proof of
Proposition~\ref{tbrwthm}. Proposition~\ref{tbrwthm} says that
%
\begin{equation}
\lim_{d(x,y)\to\infty}{\mathbf P}_u[H_y]=0\qquad \mbox{whenever
} u<\frac
{1}{2\Delta^2} \biggl(\frac{1}{\rho}-1 \biggr),
\end{equation}
and therefore, we obtain (\ref{mainstat}).
\end{pf}

\section{The connected phase}\label{sconnectedsection}
The aim of this section is to provide results regarding the phase
$u>u_c$. In Section~\ref{sucfinite} we give an example of a graph
for which this phase is nonempty. In Section~\ref{suniqmonotone}
we show that for any $u>u_c$, the interlacement set ${\mathcal I}$ is connected.

\subsection{\texorpdfstring{Finiteness of $u_c$ for $\mathbb{T}^d \times\mathbb{Z}^{d'}$}
{Finiteness of uc for Td x Zd'}}
\label{sucfinite}

In Proposition~\ref{pucfinite} below we prove that the value $u_c$ is
finite for the classical example of the product between a $d$-regular
tree and the $d'$-dimensional Euclidean lattice: $\mathbb{T}^d \times
\mathbb{Z}^{d'}$, where $d \geq3$ and $d' \geq1$. A~similar result
was proved in~\cite{GN90} for Bernoulli percolation. It is worth
mentioning that the finiteness of $u_c$ does not hold true for every
vertex-transitive nonamenable graph $G$, as we note in the following:
%
\begin{remark}\label{rtreeinf}
(1) For the infinite $d$ regular tree $\mathbb{T}^d$ ($d \geq3$), we have
%
\begin{equation}
\label{eucTd} u_c\bigl(\mathbb{T}^d\bigr) = \infty.
\end{equation}
Indeed, using Theorem 5.1, (5.7) and (5.9) of~\cite{Tei09}, we conclude
that for any $u > 0$, with ${\mathbb P}_u$-positive probability, the
root $\varnothing
\in\mathbb{T}^d$ is an isolated component of $\mathbb{T}^d \setminus
\mathcal{I}$. In this event, any two neighbors $y$ and $y'$ of
$\varnothing$ are contained in $\mathcal{I}$ ruling out
the almost sure connectedness of $\mathcal{I}$ and yielding (\ref{eucTd}).

(2) Considering again the above mentioned event, since the set
$\mathcal{I}$ has $\mathbb{P}_u$-a.s. no finite components [see
(\ref{eqW*}) and~\cite{Tei09} (2.26)] we conclude that with
positive probability there are at least two distinct infinite clusters
in $\mathcal{I}$. This together with Proposition \ref
{pnumberofclusters} gives that for
every $u > 0$ the interlacement set $\mathcal{I}$ has infinitely
many connected components ${\mathbb P}_u$-a.s.
\end{remark}

In what follows we use the same convention as in (4.1) of~\cite{Tei09}
for the product of two graphs. More precisely, if $G = (V,E)$ and $G' =
(V',E')$ are graphs, the product $G \times G'$ has vertex set $V \times
V'$, and there is an edge between $(x,x')$ and $(y,y')$ if and only if
$d(x,y) + d'(x',y') = 1$, where $d(\cdot, \cdot)$ and $d'(\cdot,
\cdot
)$ denote, respectively, the distances in $G$ and $G'$.
%
\begin{proposition}
\label{pucfinite}
For any $d \geq3$ and $d' \geq1$ we have that
%
\begin{equation}
0 < u_c \bigl( \mathbb{T}^d \times\mathbb{Z}^{d'}
\bigr) < \infty.
\end{equation}
\end{proposition}
\begin{pf}
The fact that $u_c$ is positive follows directly from Theorem \ref
{tdisconn} and the nonamenability of the graph under consideration;
see~\cite{W00}, 4.10, page 44. Therefore we focus on establishing that
$u_c < \infty$.

We first obtain a characterization of the random interlacements law on
$\mathbb{T}^d \times\mathbb{Z}^{d'}$, which resembles Theorem 5.1 in
\cite{Tei09}. For this, let us first note that the set of trajectories
\[
W^*_{\mathrm{bad}(n)} = \bigl\{w \in W^*; w \mbox{ intersects } B(\varnothing, n)
\times\mathbb{Z}^{d'} \mbox{ infinitely many times}\bigr\}
\]
has $\nu$-measure zero for all $n \geq0$, where $\varnothing$ denotes
the origin of the tree. To see why this is true, observe that the first
projection of a random walk on $\mathbb{T}^d \times\mathbb{Z}^{d'}$ is
transient. This, together with equations (\ref{eqQK}) and
(\ref{eqnuQ}), shows that $W^*_{\mathrm{bad}(n)} \cap
W^*_{B(\varnothing, n) \times B(0,m)}$ has $\nu$-measure zero for every
$m \geq0$, yielding the claim.

For $x \in V$, we write $n(x) = \inf\{n; x \in B(\varnothing,n)
\times
\mathbb{Z}^{d'}\}$ and define a map $\phi\dvtx  W^* \to\mathbb{T}^d
\times
\mathbb{Z}^{d'}$ in the following way. If $w$ belongs to some
$\{W^*_{\mathrm{bad}(n)}\}_{n \geq0}$, then $\phi(w)$ can be defined
arbitrarily, otherwise $\phi(w)$ returns the unique point $x \in
\mathbb{T}^d \times\mathbb{Z}^{d'}$ such that:
\begin{itemize}
\item$w$ does not enter $B(\varnothing, n(x)-1) \times\mathbb
{Z}^{d'}$ [trivially true if $n(x) = 1$] and
\item$x$ is the first point visited in $B(\varnothing, n(x)) \times
\mathbb{Z}^{d'}$.
\end{itemize}
Note that this is well defined out of $\bigcup_{n \geq0} W^*_{\mathrm{bad}(n)}$.

A second observation is that the sets
%
\begin{equation}
W^*_x = \bigl\{w \in W^*; \phi(w) = x \bigr\}
\end{equation}
form a disjoint partition of $W^*$. This decomposition provides us with
an alternative way to construct the set $\mathcal{I}$. More precisely,
in some auxiliar probability space $(\Omega', \mathcal{A}', P')$:
\begin{itemize}
\item we consider for every $x \in V$, independent random variables
$(J_x)_{x \in V}$, distributed as Poisson$(u \cdot\nu(W^*_x))$;
\item for $x \in V$ and $i\leq J_x$ we let $X_x^i$ and $Y_x^i$ be
random walks with distributions $P_x[\cdot|\tilde{H}_{B(\varnothing
,n(x)) \times\mathbb{Z}^{d'}} = \infty]$ and $P_x[\cdot
|H_{B(\varnothing,n(x)-1) \times\mathbb{Z}^{d'}} = \infty]$.
\end{itemize}
Using (\ref{eqQK}), (\ref{eqnuQ}) and Proposition 3.6 of~\cite{R08},
we conclude that under ${\mathbb P}_u$,
%
\begin{equation}
\label{eIurw} \mbox{the set $\mathcal{I}$ is distributed as ${ \bigcup
_{x\in V}} { \bigcup_{i \leq J_x}}
\operatorname {Range}(X^i_x) \cup
\operatorname{Range}(Y^i_x)$.}
\end{equation}

The main advantage of this representation is that the variables $J_x$
are independent, providing us with an (inhomogeneous) Bernoulli
percolation $Z_x = \mathbf{1}_{\{J_x > 0\}}$, for $x \in V$, such that
$P'[Z_x = 0] = P'[J_x = 0] = \exp\{-u \nu(W^*_x)\}$. It is also
important to note by (\ref{eqQK}) and (\ref{eqnuQ}) that
%
\begin{eqnarray}
\nu\bigl(W^*_x\bigr) & = & P_x[\tilde{H}_{B(\varnothing,n(x)) \times\mathbb
{Z}^{d'}}
= \infty, \tilde H_{\{x\}} = \infty]
d_x P_x[\tilde{H}_{B(\varnothing,n(x)-1) \times\mathbb
{Z}^{d'}} = \infty]\hspace*{-27pt}\nonumber\\[-8pt]\\[-8pt]
&\geq&\beta> 0,
\nonumber
\end{eqnarray}
uniformly on $x \in V$. This implies that $(Z_x)_{x \in V}$
stochastically dominates a Bernoulli$(1 - \exp\{-u \beta\})$ i.i.d.
percolation.

We now state a known result on the uniqueness of the infinite cluster
of Bernoulli percolation, first proved in~\cite{GN90}. More precisely,
Proposition 8.1 of~\cite{HJ06} states that for $d \geq3$, $d' \geq1$
and $p$ close enough to $1$,
%
\begin{equation}
\begin{tabular}{p{310pt}}
there is a.s. a unique infinite cluster $\mathcal{C}_\infty$
in Bernoulli($p$) site percolation on $\mathbb{T}^d \times
\mathbb {Z}^{d'}$.
\end{tabular}\hspace*{-27pt}
\end{equation}
A careful reader would observe that the result in~\cite{HJ06} is stated
for bond percolation instead. However, the second part of the proof of
Lemma 4.1 applies without modifications to site percolation, if one
uses the fact that Lemma 4.3 of~\cite{HJ06} was proved in~\cite{Sch99}
for site percolation as well. Moreover, using equation (8.9) of \cite
{HJ06} we obtain that this unique infinite cluster satisfies
%
\begin{equation}
P_x [H_{\mathcal{C}_\infty} < \infty] = 1 \qquad\mbox{for every $x \in V$}.
\end{equation}

We now take $u_o$ large enough so that for every $u \geq u_o$,
$(Z_x)_{x \in V}$ dominates a Bernoulli percolation satisfying the two
above claims. Then, for these values of $u$, the ${\mathbb P}_u$-almost
sure connectivity of $\mathcal{I}$ will follow from the
characterization in (\ref{eIurw}), proving that $u_c \leq u_o <
\infty$.
\end{pf}
%
\begin{remark}
A natural question that the above proof raises is: Is it true that for
any nonamenable graph, random interlacements dominate Bernoulli site
percolation, as we obtained in the proof of Proposition \ref
{pucfinite}? This question was posed to the authors by Itai Benjamini.

Note also that random interlacements on $\mathbb{Z}^d$ do not dominate
or become dominated by any Bernoulli site percolation; see Remark 1.1
in~\cite{SS09}.
\end{remark}

\subsection{Monotonicity of the uniqueness transition}\label{suniqmonotone}

Observe that the event $\{\mathcal{I}$ is connected$\}$ is not monotone
with respect to the set $\mathcal{I}$. This immediately raises the
question of whether there could be some $u > u_c$ for which
$\mathcal{I}^u$ is $\mathcal{P}_u$-a.s. disconnected. The next result
rules out this possibility. Note that a similar question was
considered in~\cite{HP99} in the case of Bernoulli percolation on
unimodular transitive graphs and in~\cite{Sch99} in the general
quasi-transitive case.

Fix intensities $u' > u \geq0$. For the statement of Theorem \ref
{tmonotone} below, we need to couple the random interlacement at
levels $u$ and $u'$. Actually, this is done in~\cite{Tei09} for all
values of $u \geq0$ simultaneously, but for the purpose of this paper,
it is enough to consider the following. Put ${\mathbb P}={\mathbb P}_u
\otimes{\mathbb P}_{u'-u}$. For $(\omega_1,\omega_2)\in\Omega
\times
\Omega$, put $\mathcal{I}^u={\mathcal I}(\omega_1)$ and $\mathcal
{I}^{u'-u}={\mathcal I}(\omega_2)$. Finally, we define
$\mathcal{I}^{u'}$ as $\mathcal{I}^u \cup\mathcal{I}^{u'-u}$. Then
clearly ${\mathcal I}^{u}\subset{\mathcal I}^{u'}$, and under
${\mathbb P}$, the random sets ${\mathcal I}^{u}$ and ${\mathcal
I}^{u'-u}$ are independent. Moreover, due to the Poissonian character
of random interlacements, we have that under ${\mathbb P}$, the sets
${\mathcal I}^{u}$, ${\mathcal I}^{u'}$, ${\mathcal I}^{u'-u}$ have
laws ${\mathbb P}_u$, ${\mathbb P}_{u'}$, ${\mathbb P}_{u'-u}$, respectively.
%
\begin{theorem}
\label{tmonotone}
Suppose that $G$ is such that $P_x[\tilde{H}_x = \infty]$ is uniformly
bounded from below by $\gamma> 0$. Then, for any $u' > u \geq0$, we
have that
%
\begin{equation}
\label{emonotone}\qquad \mbox{all components of $\mathcal{I}^{u'}$ contain
an infinite component of $\mathcal{I}^{u}$, $\mathbb{P}$-a.s.}
\end{equation}
\end{theorem}

Observe that transience and transitivity imply the above hypothesis.
\begin{pf*}{Proof of Theorem~\ref{tmonotone}}
To obtain the statement in (\ref{emonotone}), it is enough to prove
that, for every $x \in V$, the following event has $\mathbb
{P}$-probability one:
\begin{eqnarray*}
&&\bigl\{ \mbox{the connected component of $\mathcal{I}^{u'}$ containing
$x$}
\\
&&\hspace*{5pt}\mbox{is either empty or contains an infinite component of $
\mathcal{I}^{u}$}\bigr\}.
\end{eqnarray*}
Indeed, if this is the case, we can intersect the countable collection
of such events (as $x$ runs in $V$), and we obtain (\ref{emonotone}).

Now fix $x \in V$ and recall from (\ref{eqW*}) and (\ref
{eqinterlace}) that the set $\mathcal{I}^u$ has no finite components.
Then, since $\mathcal{I}^u \subseteq\mathcal{I}^{u'}$, the above
mentioned event equals
%
\begin{eqnarray}
\label{ehitsIu}
&&\bigl\{ %
\mbox{the connected component of $
\mathcal{I}^{u'}$ containing $x$}\nonumber\\[-8pt]\\[-8pt]
&&\hspace*{29.2pt}\mbox{is either empty or contains a
point in $\mathcal{I}^{u}$}\bigr\}.\nonumber
\end{eqnarray}

Now let us see that the above event has $\mathbb{P}$-probability one.
For this, recall that if $x \in\mathcal{I}^{u'} \setminus\mathcal
{I}^u$, then $x \in\mathcal{I}^{u'-u}$, which is equivalent to
$\omega_2(W^*_{\{x\}}) > 0$.

Given the Poissonian character of $\omega_2$, we can construct its
restriction to the set $W^*_{\{x\}}$ in the following way. First
simulate the random variable $J = \omega_2(W^*_{\{x\}})$ which has
Poisson distribution [with parameter $(u'-u)\nu(W^*_{\{x\}})$]. Then
for\vspace*{2pt} every $i = 1, \ldots, J$ we throw an independent trajectory $w \in
W^*_{\{x\}}$ with distribution\vspace*{2pt} given by $\nu$ restricted to $W^*_{\{x\}
}$ and normalized to become a probability measure. According to (\ref
{eqnuQ}), this probability distribution is given by $\pi^* \circ Q_{\{
x\}} / {e_{\{x\}}(x)}$.

Conditioned on the event $\omega_2(W^*_{\{x\}}) > 0$, we know by the
above construction and (\ref{eqQK}) that $\mathcal{I}^{u'-u}$
contains a random walk trajectory with law $P_x$. Thus, all we have to
do in order to prove that (\ref{ehitsIu}) has probability one is to
show that $\mathbb{P}$-a.s. the set $\mathcal{I}^u$ is such that
%
\begin{equation}
P_x[H_{\mathcal{I}^u} < \infty] = 1,
\end{equation}
which, by Fubini's theorem, is equivalent to proving that $P_x$-a.s.
the range of a random walk is a sequence $(x_i)_{i \geq0}$ such that
%
\begin{equation}
\label{einterinter} \mathbb{P}\bigl[ \mathcal{I}^u \cap
\{x_i\}_{i \geq0} \neq\varnothing\bigr] = 1.
\end{equation}

In fact we will exhibit a subsequence $(x_{i_j})_{j \geq0}$ of the
random walk trace, such that $\mathcal{I}^u$ almost surely intersects
$\{x_{i_j}\}_{j \geq0}$. For this, let us show that $P_x$-a.s. there
exists a subsequence $(x_{i_j})_{j \geq0}$ of the random walk trace
such that
%
\begin{equation}
\label{esubseq}
P_{x_{i_j}} [H_{\{x_{i_0}, \ldots,x_{i_{j-1}}\}} <
\infty] < \tfrac12 \qquad \mbox{for all $j \geq1$}.
\end{equation}

To see why this is true, we first claim that for every finite $K
\subset V$, the set
%
\begin{equation}
\label{eAtrans} A_{K} = \bigl\{z \in V; P_z[H_K
< \infty] \geq\tfrac12 \bigr\} \mbox{ is $P_x$-a.s. visited finitely
many times}.\hspace*{-35pt}
\end{equation}
Indeed, if the random walk had a positive probability of visiting the
set $A_{K}$ infinitely many times, we could use Borel--Cantelli's lemma
to prove that the set $K$ would also be visited infinitely often, which
is a contradiction, proving (\ref{eAtrans}). Note that for some graphs
the set $A_{K}$ could be infinite.

We can now define the subsequence $(i_j)_{j \geq0}$ that we mentioned
in (\ref{esubseq}) above. First fix a sequence $(x_i)_{i \geq0}$ that
visits $A_{B(x,n)}$ finitely often for every $n \geq0$ (note that
$A_{B(x,n)}$, $n \geq0$ is a countable family). According to (\ref
{eAtrans}) these sequences have $P_x$-probability one. Now, let $i_0 =
0$ and supposing we have defined $i_j$ for $j \leq j_o$, take $i_{j_o +
1}$ to be such that $x_{i_{j_o + 1}}$ is outside of $A_{B(x,n)}$ with
$K_{j-1}:= \{x_{i_0}, \ldots, x_{i_{j-1}}\}$ $\subset$ $B(x,n)$. This
shows (\ref{esubseq}).

We are now in position to prove (\ref{einterinter}) by considering the
following disjoint subsets of $W^*$:
%
\begin{equation}
W^*_j = \bigl\{w \in W^*; \operatorname{Range}(w) \mbox{ intersects
$x_{i_j}$ but not $K_{j-1}$}\bigr\}\qquad \mbox{for $j \geq1$}.\hspace*{-35pt}
\end{equation}
They have $\nu$-measure bounded away from $0$, as the following
calculation shows:
%
\begin{eqnarray}\quad
\nu\bigl(W^*_j\bigr) & = & \nu\bigl(W_{\{x_{i_j}\}}^* \setminus
W^*_{K_{j-1}}\bigr) \stackrel{(\rref{eqnuQ})} {=} Q_{\{x_{i_j}\}}
\bigl[(X_n)_{n \in
\mathbb{Z}} \cap K_{j-1} = \varnothing
\bigr]
\nonumber\\[-8pt]\\[-8pt]
&\stackrel{(\rref{eqQK})} {=}& P_{x_{i_j}}[H_{K_{j-1}} = \infty|
\tilde{H}_{x_{i_j}} = \infty] e_{\{x_{i_j}\}}(x_{i_j})
P_{x_{i_j}}[H_{K_{j-1}} = \infty],
\nonumber
\end{eqnarray}
which by (\ref{esubseq}) and the hypothesis of the theorem, is bounded
from below by $\gamma/ 4$.

Finally, we estimate
%
\begin{eqnarray}
\mathbb{P}\bigl[ \mathcal{I}^u \cap\{x_i
\}_{i \geq0} \neq\varnothing\bigr] &\geq& \mathbb{P}\bigl[
\mathcal{I}^u \cap\{x_{i_j}\}_{j \geq0} \neq
\varnothing\bigr]
\nonumber\\[-8pt]\\[-8pt]
&\geq& \mathbb{P}\bigl[\omega_1\bigl(W^*_j\bigr) > 0
\mbox{, for some $j \geq 1$}\bigr],
\nonumber
\end{eqnarray}
which has probability one, since the above random variables are
independent Poisson random variables with parameter bounded away from
zero. This proves (\ref{einterinter}), therefore establishing that
(\ref{ehitsIu}) has probability one, completing the proof of
Theorem~\ref{tmonotone}.
\end{pf*}

The following corollary is an immediate consequence of Theorem \ref
{tmonotone}.
%
\begin{corollary}
\label{coneuc}
For any graph $G$ satisfying the hypothesis of Theorem~\ref{tmonotone}
and for any $u > u_c$, the set $\mathcal{I}$ is $\mathbb{P}_u$-a.s.
connected. In particular, for these graphs we could have alternatively
used the definition
%
\begin{equation}
u_c = \sup\{ u \geq0; \mathcal{I} \mbox{ is not $
\mathbb{P}_u$-a.s. connected} \}.
\end{equation}
Furthermore, the connectedness transition for $\mathcal{I}$ is unique
for these graphs.
\end{corollary}
%
\begin{remark}
(1) Note that we do not necessarily suppose that the underlying graph
is transitive (or quasi-transitive) which is a standard assumption in
the case of Bernoulli percolation to obtain the monotonicity of the
uniqueness transition.

(2) In the above result we do not rule out the existence of exceptional
intensities $u > u_c$ for which $\mathcal{I}$ is disconnected.
According to Corollary~\ref{coneuc} and Fubini's theorem, the set of
such exceptional intensities must have zero Lebesgue measure.
\end{remark}

\section{Bounds for $u_*$}\label{sulowerbound}
Recall that in Theorem 4.1 of~\cite{Tei09}, it was shown that for any
nonamenable graph of bounded degree, the critical value $u_*$ is
finite. In this section, we provide a lower bound for $u_*$ which
implies the positivity of $u_*$ on any nonamenable Cayley graph. First
we recall the definition of a Cayley graph: given a finitely generated
group $H$ with symmetric generating set $S$, the (right) Cayley graph
$G=G(H,S)$ is the graph with vertex set $H$ and such that $\{u,v\}$ is
an edge if and only if $u=v s$ for some $s\in S$. Recall that the
critical value $u_c$ can be degenerated on nonamenable Cayley graphs;
see Remark~\ref{rtreeinf}. In Proposition~\ref{puthresbound} below,
we show that in contrast, $u_*$ is always nondegenerate on such graphs.
%
\begin{proposition}\label{puthresbound}
Let $G$ be a Cayley graph of degree $d$, and fix $o\in V$. Then
%
\begin{equation}
\label{euthresbound} {-\frac{1}{d P_o[\tilde{H}_o=\infty]}}\log \biggl(\frac{d}{d+\kappa_V} \biggr)
\leq u_* {\leq2 d^2 \kappa_E^{-2}\leq2
d^2 \kappa_V^{-2}.}
\end{equation}
\end{proposition}
\begin{pf}
We begin with the lower bound. Since $G$ is a Cayley graph, and the law
of ${\mathcal V}$ is invariant under $\operatorname{Aut}(G)$, Theorem 2.1 of \cite
{BLPS97} implies that if
%
\begin{equation}
\label{eblps} {\mathbb P}_u[o\in{\mathcal V}]\ge\frac{d}{d+\kappa_V},
\end{equation}
then ${\mathcal V}$ contains unbounded connected components with
positive probability. We recall that [see (\ref{evuprobx})]
%
\begin{equation}
\label{evuprob} {\mathbb P}_u[o\in{\mathcal V}]=\exp\bigl\{-u d
P_o[\tilde{H}_o=\infty]\bigr\}.
\end{equation}
We now conclude the lower bond in (\ref{euthresbound}) from (\ref
{eblps}) and (\ref{evuprob}).

We now proceed with the upper bound. The proof of Theorem 4.1 in \cite
{Tei09} shows that $u_*\le\overline{\kappa}$ where $\overline{\kappa}$ is the
constant appearing in statement (b) of Theorem 10.3 in~\cite{W00}.
However, the proof of Theorem 10.3 in~\cite{W00} shows that $\overline
{\kappa
}$ can be chosen to equal $2\kappa^2$ where $\kappa$ is the constant
appearing in the statement of Proposition~4.3 of~\cite{W00}. The
constant $\kappa$ is the same as appears in Definition 4.1 of \cite
{W00}, and comparing that definition with (\ref{ekedef}), one sees
that $\kappa=d \kappa_E^{-1}$. Thus, $u_*\le2 d^2 \kappa_E^{-2}$.
Since $|\partial_E K|\ge|\partial_V K|$ we obtain $\kappa_E\ge
\kappa_V$, and consequently $2 d^2 \kappa_E^{-2}\le2 d^2 \kappa_V^{-2}$,
completing the proof of the proposition.~%
\end{pf}

\begin{appendix}\label{app}
\section*{Appendix}

In this section we present the proof of Proposition~\ref{t01law}. We
follow here the same arguments as those of Theorem 2.1 of~\cite{Szn07},
but they are included in detail for the reader's convenience. Recall
the definitions of $t_g$, ${\mathcal Y}$ and $Q_u$ from Section \ref
{snotation}.

Before proving Proposition~\ref{t01law}, we formulate and prove
Lemma~\ref{lapproxind}. The proof of this lemma is similar to that of
Lemma 2.1 of~\cite{Bel10}, where the analogous result for random
interlacements on ${\mathbb Z}^d$ was proved; see also~\cite{Szn11},
Remark 1.5.
%
\setcounter{theorem}{0}
\begin{lemma}\label{lapproxind}
Let $u\ge0$ and $K_1$ and $K_2$ be finite disjoint subsets of
$V$. Let $F_1$ and $F_2$ be $[0,1]$-valued measurable functions on the
set of finite point-measures on $W_+$ (endowed with its canonical
$\sigma$-field).
Then
\begin{eqnarray*}
&&\bigl|{\mathbb E}_u\bigl[F_1(\mu_{K_1})
F_2(\mu_{K_2})\bigr]-{\mathbb E}_u
\bigl[F_1(\mu_{K_1})\bigr] {\mathbb E}_u
\bigl[F_2(\mu_{K_2})\bigr] \bigr|
\\
&&\qquad \le c_{u} \operatorname{cap}(K_1)\operatorname{cap}(K_2)
\mathop{\sup_{x\in K_1,}}_{y\in K_2}g(x,y).
\end{eqnarray*}
\end{lemma}
\begin{pf}
We decompose the Poisson point process $\mu_{K_1\cup K_2}$ into
four independent Poisson point processes as follows:
%
\begin{equation}
\label{emudecomp} \mu_{K_1\cup K_2}=\mu_{1,1}+\mu_{1,2}+
\mu_{2,1}+\mu_{2,2},
\end{equation}
where
\begin{eqnarray*}
\mu_{1,1}&=&\sum_{i\ge0}\delta_{w_i}1
\{X_0\in K_1,H_{K_2}=\infty\},\qquad\!\!
\mu_{1,2}=\sum_{i\ge
0}\delta_{w_i}1
\{X_0\in K_1,H_{K_2}<\infty\},
\\
\mu_{2,1}&=&\sum_{i\ge
0}\delta_{w_i}1
\{X_0\in K_2,H_{K_1}<\infty\},\qquad\!\!
\mu_{2,2}=\sum_{i\ge0}\delta_{w_i}1
\{X_0\in K_2,H_{K_1}=\infty\}.
\end{eqnarray*}
From (\ref{emukthm}) we conclude that the $\mu_{i,j}$'s are
independent Poisson point processes on~$W_+$, and their corresponding
intensity measures are given by
\begin{eqnarray*}
&&u 1\{X_0\in K_1,H_{K_2}=\infty
\}P_{e_{K_1\cup K_2}},\qquad u 1\{ X_0\in K_1,H_{K_2}<
\infty\}P_{e_{K_1\cup K_2}},
\\
&&u 1\{X_0\in K_2,H_{K_1}<\infty
\}P_{e_{K_1\cup K_2}},\qquad u 1\{X_0\in K_2,H_{K_1}=
\infty\}P_{e_{K_1\cup K_2}}.
\end{eqnarray*}
We observe that $\mu_{K_1}-\mu_{1,1}-\mu_{1,2}$ is
determined by $\mu_{2,1}$ and therefore independent of
$\mu_{1,1},\mu_{2,2}$ and $\mu_{1,2}$. In the same way,
$\mu_{K_2}-\mu_{2,2}-\mu_{2,1}$ is independent of
$\mu_{2,2},\mu_{2,1}$ and $\mu_{1,1}$. We can therefore introduce
the auxiliary Poisson processes $\mu_{2,1}'$ and
$\mu_{1,2}'$ such that they have the same law as
$\mu_{K_1}-\mu_{1,1}-\mu_{1,2}$ and
$\mu_{K_2}-\mu_{2,2}-\mu_{2,1}$, respectively, and $\mu_{2,1}'$,
$\mu_{1,2}'$, $\mu_{i,j}$, $1\le i, j\le2$ are independent.
Then
%
\begin{eqnarray}
\label{emodp1} {\mathbb E}_u\bigl[F_1(
\mu_{K_1})\bigr]&=&{\mathbb E}_u\bigl[F_1\bigl((
\mu_{K_1}-\mu_{1,1}-\mu_{1,2})+\mu_{1,1}+
\mu_{1,2}\bigr)\bigr]
\nonumber\\[-8pt]\\[-8pt]
&=&{\mathbb E}_u\bigl[F_1\bigl(\mu'_{2,1}+
\mu_{1,1}+\mu_{1,2}\bigr)\bigr]\nonumber
\end{eqnarray}
and in the same way,
%
\begin{equation}
\label{emodp2} {\mathbb E}_u\bigl[F_2(
\mu_{K_2})\bigr]={\mathbb E}_u\bigl[F_2\bigl(
\mu'_{1,2}+\mu_{2,2}+\mu_{2,1}\bigr)
\bigr].
\end{equation}
Using (\ref{emodp1}), (\ref{emodp2}) and the independence of the
Poisson processes $\mu'_{2,1}+\mu_{1,1}+\mu_{1,2}$ and
$\mu'_{1,2}+\mu_{2,2}+\mu_{2,1}$, we get
%
\begin{eqnarray}
\label{eprod}
&&{\mathbb E}_u\bigl[F_1(
\mu_{K_1})\bigr] {\mathbb E}_u\bigl[F_2(
\mu_{K_2})\bigr]\nonumber\\[-8pt]\\[-8pt]
&&\qquad={\mathbb E}_u\bigl[F_1\bigl(
\mu'_{2,1}+\mu_{1,1}+\mu_{1,2}
\bigr)F_2\bigl(\mu'_{1,2}+\mu_{2,2}+
\mu_{2,1}\bigr)\bigr].\nonumber
\end{eqnarray}
From (\ref{eprod}) we see that
%
\begin{eqnarray}
\label{edepend}\qquad
&&\bigl|{\mathbb E}_u\bigl[F_1(
\mu_{K_1}) F_2(\mu_{K_2})\bigr]-{\mathbb
E}_u\bigl[F_1(\mu_{K_1})\bigr] {\mathbb
E}_u\bigl[F_2(\mu_{K_2})\bigr]\bigr|
\nonumber
\\
&&\qquad\le {\mathbb P}_u\bigl[\mu'_{2,1}\neq0
\mbox{ or }\mu'_{1,2}\neq0\mbox{ or }\mu_{2,1}
\neq0\mbox{ or }\mu_{1,2}\neq0\bigr]
\nonumber\\[-8pt]\\[-8pt]
&&\qquad\le 2\bigl({\mathbb P}_u[\mu_{2,1}\neq0]+{\mathbb
P}_u[\mu_{1,2}\neq0]\bigr)
\nonumber
\\
&&\qquad\le 2u\bigl(P_{e_{K_1\cup K_2}}[X_0\in K_1,H_{K_2}<
\infty] +P_{e_{K_1\cup K_2}}[X_0\in K_2,H_{K_1}<
\infty]\bigr).\nonumber
\end{eqnarray}
We now bound the two last terms in the above equation,
%
\begin{eqnarray}
P_{e_{K_1\cup K_2}}[X_0\in K_1,H_{K_2}<\infty]
&\le&\sum_{x\in K_1} e_{K_1}(x)P_x[H_{K_2}<
\infty]
\nonumber
\\
&=&\sum_{x\in K_1,y\in K_2} e_{K_1}(x) g(x,y)
e_{K_2}(y)
\\
&\le& \operatorname{cap}(K_1) \operatorname{cap}(K_2)
\sup_{x\in K_1,
y\in
K_2} g(x,y).
\nonumber
\end{eqnarray}
A similar estimate holds for $P_{e_{K_1\cup K_2}}[X_0\in
K_2,H_{K_1}<\infty]$, and the result follows.
\end{pf}
\begin{pf*}{Proof of Proposition~\ref{t01law}}
We follow the proof of Theorem 2.1 in~\cite{Szn07}, which goes through
with only minor modifications. We define the map $\psi\dvtx  \Omega\to\{
0,1\}^V$ given by $\psi(\omega)=(1\{x\in{\mathcal V}(\omega)\}
)_{x\in
V}$. Then $Q_u=\psi\circ{\mathbb P}_u$ and moreover
%
\begin{equation}
\label{eshifteq} t_g\circ\psi=\psi\circ\tau_g,\qquad g\in\autg.
\end{equation}

Choose a sequence of vertices $v,v_1,v_2,\ldots\in V$ such that
$d(v,v_i)\to\infty$ as $i\to\infty$. For each $i\ge1$, let $g_i\in
\autg$ be such that $g_i(v)=v_i$. These choices are possible due to our
assumption that $G$ is an infinite transitive graph. To prove the
ergodicity statement, it suffices to establish that for any finite
$K\subset V$ and any $[0,1]$-valued $\sigma(Y_z,z\in K)$-measurable
function $f$ on $\{0,1\}^{{\mathbb Z}^d}$, the following limit holds:
%
\begin{equation}
\label{eexplimit} \lim_{i\to\infty}E^{Q_u}[f f\circ
\tau_{g_i}]=E^{Q_u}[f]^2.
\end{equation}

Bound (\ref{eexplimit}) gives the mixing of the interlacements set
with respect to the automorphisms of $G$. The $\{0,1\}$-law (\ref
{e01law}) can be classically deduced as follows: if $A\in{\mathcal Y}$
is invariant under $\operatorname{Aut}(G)$, then one can do $L^1(Q_u)$-approximation
of its indicator function by functions $f$ as above. With (\ref
{eexplimit}), one obtains in a standard way that $Q_u(A)=Q_u(A)^2$, so
that $Q_u(A)\in\{0,1\}$.

Using (\ref{eshifteq}), equation (\ref{eexplimit}) will follow if we
show that for any finite $K\subset V$ and any $[0,1]$-valued function
$F$ on the set of finite point-measures on $W_+$,
%
\begin{equation}
\label{eexplimit2} \lim_{i\to\infty}{\mathbb E}_u\bigl[F(
\mu_K) F(\mu_K)\circ\tau_{g_i}\bigr]={\mathbb
E}_u\bigl[F(\mu_K)\bigr]^2.
\end{equation}
However, ${\mathbb E}_u[F(\mu_K) F(\mu_K)\circ\tau_{g_i}]={\mathbb
E}_u[F(\mu_K) H(\mu_{g_i(K)})]$ for some $H$ in the same class of
functions as $F$. Since $d(K,g_i(K))\to\infty$ as $i\to\infty$ and $G$
is transient, we conclude that $\sup_{x\in K, y\in g_i(K)}$
$g(x,y)\to
0$ as $i\to\infty$. An appeal to Lemma~\ref{lapproxind} now gives
(\ref{eexplimit2}), and the theorem follows.
\end{pf*}
\end{appendix}

\section*{Acknowledgment}

The authors thank Itai Benjamini for suggesting some of the problems
dealt with in the paper.



\printaddresses

\end{document}